\documentclass[12pt]{article}
\usepackage{amsmath,amssymb,a4wide,color,url,comment}
\usepackage{wrapfig}
\usepackage{scalerel,stackengine}
\stackMath
\renewcommand\widehat[1]{%
\savestack{\tmpbox}{\stretchto{%
  \scaleto{%
    \scalerel*[\widthof{\ensuremath{#1}}]{\kern-.6pt\bigwedge\kern-.6pt}%
    {\rule[-\textheight/2]{1ex}{\textheight}}
  }{\textheight}%
}{0.5ex}}%
\stackon[1pt]{#1}{\tmpbox}%
}\definecolor{grey}{rgb}{0.7, 0.7, 0.7}
\newtheorem{lemma}{Lemma}[section]
\newtheorem{theorem}{Theorem}[section]

\newcommand\Ltwow{{L^2_{c^2/c_0^2}(\Omega)}}

\newcommand\ssl{\mathfrak{s}}
\newcommand\nlc{\eta}
\newcommand\kappaintro{\kappa}

\newcommand\dssl{\underline{d\ssl}}
\newcommand\dnlc{\underline{d\nlc}}
\newcommand\du{\underline{du}}
\newcommand\delssl{\underline{\delta\ssl}}
\newcommand\delnlc{\underline{\delta\nlc}}
\newcommand\delu{\underline{\delta u}}

\newcommand\Lapl[1]{\widehat{#1}}
\newcommand\regpar{\gamma}
\def\FZ{{\sc fz\,}}
\def\CWKCH{{\sc ch\,}}

\newcommand{\Margin}[1]{}
\newcommand{\anapap}[1]{}
\newcommand{\revision}[1]{#1}
\newcommand{\arxiv}[1]{#1}
\newcommand{\rerevision}[1]{#1}

\newbox\figurelegendone
\newbox\figurelegendtwo
\newbox\figurelegendthree
\newbox\figurelegendfour
\newbox\figureone
\newbox\figuretwo
\newbox\figurethree 
\newbox\figurefour
\newbox\figurefive
\newbox\figuresix
\newbox\figureseven
\newbox\figureeight
\newbox\figurenine

\input colordvi
\font\tenrm=cmr10
\font\teni=cmmi10 \skewchar\teni='177
\font\tensy=cmsy10 \skewchar\tensy='60
\font\tenex=cmex10
\font\tenit=cmti10
\font\tensl=cmsl10
\font\tenbf=cmbx10
\font\tentt=cmtt10
\font\ninerm=cmr9
\font\ninei=cmmi9 \skewchar\ninei='177
\font\ninesy=cmsy9 \skewchar\ninesy='60
\font\nineit=cmti9
\font\ninesl=cmsl9
\font\ninebf=cmbx9
\font\ninett=cmtt9
\font\eightrm=cmr8
\font\eighti=cmmi8 \skewchar\eighti='177
\font\eightsy=cmsy8 \skewchar\eightsy='60
\font\eightit=cmti8
\font\eightsl=cmsl8
\font\eightbf=cmbx8
\font\eighttt=cmtt8
\font\sevenrm=cmr7
\font\seveni=cmmi7 \skewchar\seveni='177
\font\sevensy=cmsy7 \skewchar\sevensy='60
\font\sevenbf=cmbx7
\font\sevenit=cmmi7
\font\sevensl=cmmi7
\font\seventt=cmr7
\font\sixrm=cmr6
\font\sixi=cmmi6 \skewchar\sixi='177
\font\sixsy=cmsy6 \skewchar\sixsy='60
\font\sixbf=cmbx6
\font\fiverm=cmr5
\font\fivei=cmmi5 \skewchar\fivei='177
\font\fivesy=cmsy5 \skewchar\fivesy='60
\font\fivebf=cmbx5
\def\tenpoint{\def\rm{\fam0\tenrm}%
        \textfont0=\tenrm \scriptfont0=\sevenrm \scriptscriptfont0=\fiverm
        \textfont1=\teni \scriptfont1=\seveni \scriptscriptfont1=\fivei
        \textfont2=\tensy \scriptfont2=\sevensy \scriptscriptfont2=\fivesy
        \textfont3=\tenex \scriptfont3=\tenex \scriptscriptfont3=\tenex
        \def\it{\fam\itfam\tenit}%
        \textfont\itfam=\tenit
        \def\sl{\fam\slfam\tensl}%
        \textfont\slfam=\tensl
        \def\bf{\fam\bffam\tenbf}%
        \textfont\bffam=\tenbf \scriptfont\bffam=\sevenbf
                \scriptscriptfont\bffam=\fivebf
        \def\tt{\fam\ttfam\tentt}%
        \textfont\ttfam=\tentt
        \normalbaselineskip=12pt%
        \let\sc=\eightrm        
        \setbox\strutbox=\hbox{\vrule height8.5pt depth3.5pt width0pt}%
        \normalbaselines\rm}
\def\ninepoint{\def\rm{\fam0\ninerm}%
        \textfont0=\ninerm \scriptfont0=\sixrm \scriptscriptfont0=\fiverm
        \textfont1=\ninei \scriptfont1=\sixi \scriptscriptfont1=\fivei
        \textfont2=\ninesy \scriptfont2=\sixsy \scriptscriptfont2=\fivesy
        \textfont3=\tenex \scriptfont3=\tenex \scriptscriptfont3=\tenex
        \def\it{\fam\itfam\nineit}%
        \textfont\itfam=\nineit
        \def\sl{\fam\slfam\ninesl}%
        \textfont\slfam=\ninesl
        \def\bf{\fam\bffam\ninebf}%
        \textfont\bffam=\ninebf \scriptfont\bffam=\sixbf
                \scriptscriptfont\bffam=\fivebf
        \def\tt{\fam\ttfam\ninett}%
        \textfont\ttfam=\ninett
        \normalbaselineskip=11pt%
        \let\sc=\sevenrm        
        \setbox\strutbox=\hbox{\vrule height8pt depth3pt width0pt}%
        \normalbaselines\rm}
\def\eightpoint{\def\rm{\fam0\eightrm}%
        \textfont0=\eightrm \scriptfont0=\sixrm \scriptscriptfont0=\fiverm
        \textfont1=\eighti \scriptfont1=\sixi \scriptscriptfont1=\fivei
        \textfont2=\eightsy \scriptfont2=\sixsy \scriptscriptfont2=\fivesy
        \textfont3=\tenex \scriptfont3=\tenex \scriptscriptfont3=\tenex
        \def\it{\fam\itfam\eightit}%
        \textfont\itfam=\eightit
        \def\sl{\fam\slfam\eightsl}%
        \textfont\slfam=\eightsl
        \def\bf{\fam\bffam\eightbf}%
        \textfont\bffam=\eightbf \scriptfont\bffam=\sixbf
                \scriptscriptfont\bffam=\fivebf
        \def\tt{\fam\ttfam\eighttt}%
        \textfont\ttfam=\eighttt
        \normalbaselineskip=9pt%
        \let\sc=\sixrm  
        \setbox\strutbox=\hbox{\vrule height7pt depth2pt width0pt}%
        \normalbaselines\rm}
\def\sevenpoint{\def\rm{\fam0\sevenrm}%
        \textfont0=\sevenrm \scriptfont0=\fiverm \scriptscriptfont0=\fiverm
        \textfont1=\seveni \scriptfont1=\fivei \scriptscriptfont1=\fivei
        \textfont2=\sevensy \scriptfont2=\fivesy \scriptscriptfont2=\fivesy
        \textfont3=\tenex \scriptfont3=\tenex \scriptscriptfont3=\tenex
        \def\it{\fam\itfam\sevenit}%
        \textfont\itfam=\sevenit
        \def\sl{\fam\slfam\sevensl}%
        \textfont\slfam=\sevensl
        \def\bf{\fam\bffam\sevenbf}%
        \textfont\bffam=\sevenbf \scriptfont\bffam=\fivebf
                \scriptscriptfont\bffam=\fivebf
        \def\tt{\fam\ttfam\seventt}%
        \textfont\ttfam=\seventt
        \normalbaselineskip=8pt%
        \let\sc=\fiverm  
        \setbox\strutbox=\hbox{\vrule height6pt depth2pt width0pt}%
        \normalbaselines\rm}

\input pictex

\newdimen\xfiglen \newdimen\yfiglen
\setbox\figurelegendone=\hbox{
\beginpicture
  \setcoordinatesystem units <0.08\xfiglen,0.014\yfiglen> 
  \setplotarea x from 0 to 2, y from 0 to 3
\linethickness=0.6pt
\scriptsize
   \setsolid  
   \putrule from 0 3 to 1 3 
   \setdashes <3pt> 
   \putrule from 0 2 to 1 2 
   \setplotsymbol ({{\tiny $\ast$}}) \plotsymbolspacing=4pt
   \plot 
      0 1
      1 1
   /
   \setplotsymbol ({{\tiny $\bullet$}}) \plotsymbolspacing=4pt
   \plot 
      0 0
      1 0
   /
   \setdots <3pt> 
   \setdots
   \putrule from 0 0 to 1 0
  \put {actual}  [l] at 1.2 3
  \put {iter 1}  [l] at 1.2 2
  \put {iter 2}  [l] at 1.2 1
  \put {iter 3}  [l] at 1.2 0
\endpicture
}
\setbox\figurelegendtwo=\hbox{
\beginpicture
  \setcoordinatesystem units <0.08\xfiglen,0.014\yfiglen> 
  \setplotarea x from 0 to 2, y from 0 to 4
\linethickness=0.6pt
\scriptsize
   \setsolid  
   \putrule from 0 4 to 1 4 
   \setdashes <3pt> 
   \putrule from 0 3 to 1 3 
   \setplotsymbol ({{\tiny $\ast$}}) \plotsymbolspacing=4pt
   \plot 
      0 2
      1 2
   /
   \setplotsymbol ({{\tiny $\circ$}}) \plotsymbolspacing=4pt
   \plot 
      0 1
      1 1
   /
   \setplotsymbol ({{\tiny $\bullet$}}) \plotsymbolspacing=4pt
\putrule from 0 0 to 1 0
   \plot 
      0 0
      1 0
   /
   \setdots <3pt> 
   \setdots
   \putrule from 0 0 to 1 0
\scriptsize
  \put {actual}  [l] at 1.2 4
  \put {iter 1}  [l] at 1.2 3
  \put {iter 2}  [l] at 1.2 2
  \put {iter 5}  [l] at 1.2 1
  \put {iter 8}  [l] at 1.2 0
\endpicture
}
\setbox\figurelegendthree=\hbox{
\beginpicture
  \setcoordinatesystem units <0.08\xfiglen,0.014\yfiglen> 
  \setplotarea x from 0 to 2, y from 0 to 4
\linethickness=0.6pt
\scriptsize
   \setsolid  
   \putrule from 0 4 to 1 4 
   \setdashes <3pt> 
   \putrule from 0 3 to 1 3 
   \setplotsymbol ({{\tiny $\ast$}}) \plotsymbolspacing=4pt
   \plot 
      0 2
      1 2
   /
   \setplotsymbol ({{\tiny $\circ$}}) \plotsymbolspacing=4pt
   \plot 
      0 1
      1 1
   /
   \setplotsymbol ({{\tiny $\bullet$}}) \plotsymbolspacing=4pt
\putrule from 0 0 to 1 0
   \plot 
      0 0
      1 0
   /
   \setdots <3pt> 
   \setdots
   \putrule from 0 0 to 1 0
\scriptsize
  \put {actual}  [l] at 1.2 4
  \put {iter 1}  [l] at 1.2 3
  \put {iter 2}  [l] at 1.2 2
  \put {iter 3}  [l] at 1.2 1
  \put {iter 4}  [l] at 1.2 0
\endpicture
}
\setbox\figurelegendfour=\hbox{
\beginpicture
  \setcoordinatesystem units <0.08\xfiglen,0.014\yfiglen> 
  \setplotarea x from 0 to 2, y from 0 to 4
\linethickness=0.6pt
\scriptsize
   \setsolid  
   \putrule from 0 4 to 1 4 
   \setdashes <3pt> 
   \putrule from 0 3 to 1 3 
   \setplotsymbol ({{\tiny $\ast$}}) \plotsymbolspacing=4pt
   \plot 
      0 2
      1 2
   /
   \setplotsymbol ({{\tiny $\circ$}}) \plotsymbolspacing=4pt
   \plot 
      0 1
      1 1
   /
   \setplotsymbol ({{\tiny $\bullet$}}) \plotsymbolspacing=4pt
\putrule from 0 0 to 1 0
   \plot 
      0 0
      1 0
   /
   \setdots <3pt> 
   \setdots
   \putrule from 0 0 to 1 0
\scriptsize
  \put {actual}  [l] at 1.2 4
  \put {iter 2}  [l] at 1.2 3
  \put {iter 4}  [l] at 1.2 2
  \put {iter 6}  [l] at 1.2 1
  \put {iter 8}  [l] at 1.2 0
\endpicture
}
\setbox\figureone=\vbox{\hsize=\xfiglen  
\beginpicture
  \setcoordinatesystem units <\xfiglen,\yfiglen>  point at 0.0 0.0
  \setplotarea x from 0 to 1, y from 0 to 0.2
\scriptsize
  \axis bottom shiftedto y=0.0 ticks short numbered from 0 to 1 by 0.25 /
  \axis left ticks short numbered from 0.0 to 0.2 by 0.1 /
\footnotesize
\put {$\eta(x)$} [l] at 0.03 0.2
\linethickness=0.9pt
\setlinear
\setsolid
\relax
\plot  
    0.00    0
    0.05    0
    0.10    0
    0.15    0
    0.20    0
    0.25    0
    0.30    0
    0.40    0.030
    0.50    0.045
    0.60    0.075
    0.70    0.150
    0.80    0.150
    0.90    0.100
    0.95    0.030
    1.00    0 /
\linethickness=0.3pt
\setdashes <3pt>
\setquadratic
\relax
\plot 0         0
    0.0200    0.0136
    0.0400    0.0272
    0.0600    0.0406
    0.0800    0.0541
    0.1000    0.0672
    0.1200    0.0798
    0.1400    0.0922
    0.1600    0.1034
    0.1800    0.1146
    0.2000    0.1237
    0.2200    0.1326
    0.2400    0.1396
    0.2600    0.1454
    0.2800    0.1499
    0.3000    0.1520
    0.3200    0.1537
    0.3400    0.1519
    0.3600    0.1501
    0.3800    0.1457
    0.4000    0.1406
    0.4200    0.1345
    0.4400    0.1276
    0.4600    0.1205
    0.4800    0.1132
    0.5000    0.1059
    0.5200    0.0996
    0.5400    0.0934
    0.5600    0.0889
    0.5800    0.0851
    0.6000    0.0827
    0.6200    0.0821
    0.6400    0.0818
    0.6600    0.0833
    0.6800    0.0848
    0.7000    0.0866
    0.7200    0.0885
    0.7400    0.0896
    0.7600    0.0902
    0.7800    0.0899
    0.8000    0.0878
    0.8200    0.0855
    0.8400    0.0801
    0.8600    0.0748
    0.8800    0.0668
    0.9000    0.0583
    0.9200    0.0483
    0.9400    0.0373
    0.9600    0.0257
    0.9800    0.0129
    1.0000         0 /

\setplotsymbol ({{\tiny $\ast$}}) \plotsymbolspacing=8pt
\plot 0         0
    0.0200         0
    0.0400         0
    0.0600         0
    0.0800         0
    0.1000         0
    0.1200         0
    0.1400         0
    0.1600         0
    0.1800         0
    0.2000         0
    0.2200         0
    0.2400         0
    0.2600         0
    0.2800         0
    0.3000         0
    0.3200         0
    0.3400         0
    0.3600         0
    0.3800         0
    0.4000         0
    0.4200         0
    0.4400         0
    0.4600         0
    0.4800    0.0004
    0.5000    0.0048
    0.5200    0.0114
    0.5400    0.0180
    0.5600    0.0266
    0.5800    0.0359
    0.6000    0.0465
    0.6200    0.0586
    0.6400    0.0709
    0.6600    0.0840
    0.6800    0.0972
    0.7000    0.1095
    0.7200    0.1216
    0.7400    0.1323
    0.7600    0.1421
    0.7800    0.1505
    0.8000    0.1563
    0.8200    0.1616
    0.8400    0.1616
    0.8600    0.1616
    0.8800    0.1543
    0.9000    0.1453
    0.9200    0.1291
    0.9400    0.1069
    0.9600    0.0798
    0.9800    0.0399
    1.0000         0 /

\setplotsymbol ({{\tiny $\bullet$}}) \plotsymbolspacing=8pt
\setdots
\plot 0         0
    0.0200         0
    0.0400         0
    0.0600         0
    0.0800         0
    0.1000         0
    0.1200         0
    0.1400         0
    0.1600         0
    0.1800         0
    0.2000         0
    0.2200         0
    0.2400         0
    0.2600         0
    0.2800         0
    0.3000         0
    0.3200         0
    0.3400         0
    0.3600    0.0051
    0.3800    0.0153
    0.4000    0.0255
    0.4200    0.0346
    0.4400    0.0429
    0.4600    0.0505
    0.4800    0.0567
    0.5000    0.0670
    0.5200    0.0770
    0.5400    0.0869
    0.5600    0.0972
    0.5800    0.1076
    0.6000    0.1181
    0.6200    0.1288
    0.6400    0.1392
    0.6600    0.1477
    0.6800    0.1562
    0.7000    0.1600
    0.7200    0.1632
    0.7400    0.1626
    0.7600    0.1597
    0.7800    0.1546
    0.8000    0.1454
    0.8200    0.1359
    0.8400    0.1225
    0.8600    0.1091
    0.8800    0.0941
    0.9000    0.0788
    0.9200    0.0633
    0.9400    0.0476
    0.9600    0.0318
    0.9800    0.0159
    1.0000         0 /

\endpicture
} 
\setbox\figuretwo=\vbox{\hsize=\xfiglen  
\beginpicture
  \setcoordinatesystem units <\xfiglen,\yfiglen>  point at 0.0 0.0
  \setplotarea x from 0 to 1, y from 0 to 0.2
\scriptsize
  \axis bottom shiftedto y=0.0 ticks short numbered from 0 to 1 by 0.25 /
  \axis left ticks short numbered from 0.0 to 0.2 by 0.1 /
\footnotesize
\put {\copy\figurelegendone} [rt] at 1.05 0.2
\put {$s(x)-1$} [l] at 0.03 0.2
\setsolid
\linethickness=0.9pt
\setlinear
\plot
         0         0
    0.2000         0
    0.4000    0.1100
    0.6000    0.1300
    0.8000    0.0400
    1.0000         0
/

\linethickness=0.3pt
\setdashes <3pt>
\plot
         0    0.0000
    0.0200    0.0015
    0.0400    0.0031
    0.0600    0.0045
    0.0800    0.0059
    0.1000    0.0080
    0.1200    0.0111
    0.1400    0.0147
    0.1600    0.0212
    0.1800    0.0277
    0.2000    0.0366
    0.2200    0.0458
    0.2400    0.0561
    0.2600    0.0672
    0.2800    0.0791
    0.3000    0.0922
    0.3200    0.1055
    0.3400    0.1202
    0.3600    0.1349
    0.3800    0.1508
    0.4000    0.1668
    0.4200    0.1822
    0.4400    0.1969
    0.4600    0.2103
    0.4800    0.2200
    0.5000    0.2298
    0.5200    0.2339
    0.5400    0.2381
    0.5600    0.2383
    0.5800    0.2370
    0.6000    0.2327
    0.6200    0.2247
    0.6400    0.2154
    0.6600    0.2000
    0.6800    0.1846
    0.7000    0.1670
    0.7200    0.1491
    0.7400    0.1316
    0.7600    0.1143
    0.7800    0.0977
    0.8000    0.0826
    0.8200    0.0677
    0.8400    0.0559
    0.8600    0.0441
    0.8800    0.0351
    0.9000    0.0267
    0.9200    0.0194
    0.9400    0.0130
    0.9600    0.0073
    0.9800    0.0037
    1.0000    0.0000
/ 

\setplotsymbol ({{\tiny $\ast$}}) \plotsymbolspacing=8pt
\plot
         0    0.0000
    0.0200    0.0145
    0.0400    0.0290
    0.0600    0.0419
    0.0800    0.0541
    0.1000    0.0653
    0.1200    0.0753
    0.1400    0.0851
    0.1600    0.0940
    0.1800    0.1028
    0.2000    0.1077
    0.2200    0.1123
    0.2400    0.1137
    0.2600    0.1133
    0.2800    0.1122
    0.3000    0.1096
    0.3200    0.1072
    0.3400    0.1053
    0.3600    0.1034
    0.3800    0.1052
    0.4000    0.1079
    0.4200    0.1119
    0.4400    0.1171
    0.4600    0.1215
    0.4800    0.1237
    0.5000    0.1260
    0.5200    0.1246
    0.5400    0.1233
    0.5600    0.1198
    0.5800    0.1155
    0.6000    0.1090
    0.6200    0.0997
    0.6400    0.0892
    0.6600    0.0730
    0.6800    0.0568
    0.7000    0.0414
    0.7200    0.0260
    0.7400    0.0135
    0.7600    0.0026
    0.7800    0.0000
    0.8000    0.0000
    0.8200    0.0000
    0.8400    0.0000
    0.8600    0.0000
    0.8800    0.0000
    0.9000    0.0000
    0.9200    0.0000
    0.9400    0.0000
    0.9600    0.0000
    0.9800    0.0000
    1.0000    0.0000
/

\setplotsymbol ({{\tiny $\bullet$}}) \plotsymbolspacing=8pt
\setdots
\plot
         0    0.0000
    0.0200    0.0000
    0.0400    0.0000
    0.0600    0.0000
    0.0800    0.0000
    0.1000    0.0000
    0.1200    0.0000
    0.1400    0.0000
    0.1600    0.0000
    0.1800    0.0000
    0.2000    0.0006
    0.2200    0.0071
    0.2400    0.0118
    0.2600    0.0153
    0.2800    0.0182
    0.3000    0.0198
    0.3200    0.0214
    0.3400    0.0231
    0.3600    0.0249
    0.3800    0.0315
    0.4000    0.0392
    0.4200    0.0492
    0.4400    0.0610
    0.4600    0.0721
    0.4800    0.0811
    0.5000    0.0902
    0.5200    0.0963
    0.5400    0.1025
    0.5600    0.1075
    0.5800    0.1122
    0.6000    0.1145
    0.6200    0.1141
    0.6400    0.1120
    0.6600    0.1018
    0.6800    0.0916
    0.7000    0.0805
    0.7200    0.0693
    0.7400    0.0590
    0.7600    0.0492
    0.7800    0.0464
    0.8000    0.0439
    0.8200    0.0412
    0.8400    0.0368
    0.8600    0.0323
    0.8800    0.0276
    0.9000    0.0227
    0.9200    0.0177
    0.9400    0.0127
    0.9600    0.0079
    0.9800    0.0039
/

\endpicture
} 
\setbox\figurethree=\vbox{\hsize=\xfiglen  
\beginpicture
  \setcoordinatesystem units <\xfiglen,\yfiglen>  point at 0.0 0.0
  \setplotarea x from 0 to 1, y from 0 to 0.2
\scriptsize
  \axis bottom shiftedto y=0.0 ticks short numbered from 0 to 1 by 0.25 /
  \axis left ticks short numbered from 0.0 to 0.2 by 0.1 /
\footnotesize
\put {$\eta(x)$} [l] at 0.03 0.2
\linethickness=0.9pt
\setlinear
\setsolid
\relax
\plot
    0.00    0
    0.05    0
    0.10    0
    0.15    0
    0.20    0
    0.25    0
    0.30    0
    0.40    0.030
    0.50    0.045
    0.60    0.075
    0.70    0.150
    0.80    0.150
    0.90    0.100
    0.95    0.030
    1.00    0 /
\linethickness=0.3pt
\setdashes <3pt>
\setquadratic
\relax
\plot
         0         0
    0.0200    0.0137
    0.0400    0.0273
    0.0600    0.0409
    0.0800    0.0544
    0.1000    0.0676
    0.1200    0.0804
    0.1400    0.0928
    0.1600    0.1039
    0.1800    0.1150
    0.2000    0.1239
    0.2200    0.1326
    0.2400    0.1393
    0.2600    0.1448
    0.2800    0.1491
    0.3000    0.1511
    0.3200    0.1527
    0.3400    0.1512
    0.3600    0.1496
    0.3800    0.1455
    0.4000    0.1409
    0.4200    0.1352
    0.4400    0.1285
    0.4600    0.1218
    0.4800    0.1145
    0.5000    0.1073
    0.5200    0.1012
    0.5400    0.0950
    0.5600    0.0906
    0.5800    0.0868
    0.6000    0.0843
    0.6200    0.0833
    0.6400    0.0828
    0.6600    0.0839
    0.6800    0.0851
    0.7000    0.0867
    0.7200    0.0884
    0.7400    0.0892
    0.7600    0.0894
    0.7800    0.0887
    0.8000    0.0860
    0.8200    0.0831
    0.8400    0.0769
    0.8600    0.0707
    0.8800    0.0620
    0.9000    0.0528
    0.9200    0.0428
    0.9400    0.0320
    0.9600    0.0212
    0.9800    0.0106
    1.0000         0
/

\setplotsymbol ({{\tiny $\ast$}}) \plotsymbolspacing=8pt
\setdots
\plot
         0         0
    0.0200         0
    0.0400         0
    0.0600         0
    0.0800         0
    0.1000         0
    0.1200         0
    0.1400         0
    0.1600         0
    0.1800         0
    0.2000         0
    0.2200         0
    0.2400         0
    0.2600         0
    0.2800         0
    0.3000         0
    0.3200         0
    0.3400         0
    0.3600         0
    0.3800         0
    0.4000         0
    0.4200         0
    0.4400         0
    0.4600         0
    0.4800    0.0025
    0.5000    0.0081
    0.5200    0.0152
    0.5400    0.0223
    0.5600    0.0309
    0.5800    0.0400
    0.6000    0.0500
    0.6200    0.0609
    0.6400    0.0721
    0.6600    0.0839
    0.6800    0.0958
    0.7000    0.1072
    0.7200    0.1186
    0.7400    0.1287
    0.7600    0.1382
    0.7800    0.1462
    0.8000    0.1519
    0.8200    0.1570
    0.8400    0.1568
    0.8600    0.1565
    0.8800    0.1489
    0.9000    0.1396
    0.9200    0.1235
    0.9400    0.1017
    0.9600    0.0755
    0.9800    0.0378
    1.0000         0
/

\setplotsymbol ({{\tiny $\circ$}}) \plotsymbolspacing=8pt
\setdots
\plot
         0         0
    0.0200         0
    0.0400         0
    0.0600         0
    0.0800         0
    0.1000         0
    0.1200         0
    0.1400         0
    0.1600         0
    0.1800         0
    0.2000         0
    0.2200         0
    0.2400         0
    0.2600         0
    0.2800         0
    0.3000         0
    0.3200         0
    0.3400         0
    0.3600         0
    0.3800         0
    0.4000    0.0008
    0.4200    0.0096
    0.4400    0.0167
    0.4600    0.0229
    0.4800    0.0290
    0.5000    0.0381
    0.5200    0.0467
    0.5400    0.0554
    0.5600    0.0651
    0.5800    0.0752
    0.6000    0.0862
    0.6200    0.0981
    0.6400    0.1101
    0.6600    0.1220
    0.6800    0.1339
    0.7000    0.1431
    0.7200    0.1520
    0.7400    0.1572
    0.7600    0.1602
    0.7800    0.1605
    0.8000    0.1559
    0.8200    0.1507
    0.8400    0.1388
    0.8600    0.1269
    0.8800    0.1107
    0.9000    0.0937
    0.9200    0.0755
    0.9400    0.0564
    0.9600    0.0374
    0.9800    0.0187
    1.0000         0
/

\setplotsymbol ({{\tiny $\bullet$}}) \plotsymbolspacing=8pt
\setdots
\plot
         0         0
    0.0200         0
    0.0400         0
    0.0600         0
    0.0800         0
    0.1000         0
    0.1200         0
    0.1400         0
    0.1600         0
    0.1800         0
    0.2000         0
    0.2200         0
    0.2400         0
    0.2600         0
    0.2800         0
    0.3000         0
    0.3200    0.0038
    0.3400    0.0089
    0.3600    0.0140
    0.3800    0.0171
    0.4000    0.0206
    0.4200    0.0303
    0.4400    0.0370
    0.4600    0.0421
    0.4800    0.0454
    0.5000    0.0517
    0.5200    0.0567
    0.5400    0.0618
    0.5600    0.0684
    0.5800    0.0755
    0.6000    0.0842
    0.6200    0.0948
    0.6400    0.1055
    0.6600    0.1173
    0.6800    0.1291
    0.7000    0.1386
    0.7200    0.1480
    0.7400    0.1534
    0.7600    0.1566
    0.7800    0.1568
    0.8000    0.1515
    0.8200    0.1455
    0.8400    0.1322
    0.8600    0.1188
    0.8800    0.1015
    0.9000    0.0833
    0.9200    0.0650
    0.9400    0.0466
    0.9600    0.0293
    0.9800    0.0146
    1.0000         0
/

\endpicture
} 
\setbox\figurefour=\vbox{\hsize=\xfiglen  
\beginpicture
  \setcoordinatesystem units <\xfiglen,\yfiglen>  point at 0.0 0.0
  \setplotarea x from 0 to 1, y from 0 to 0.2
\scriptsize
  \axis bottom shiftedto y=0.0 ticks short numbered from 0 to 1 by 0.25 /
  \axis left ticks short numbered from 0.0 to 0.2 by 0.1 /
\footnotesize
\put {\copy\figurelegendtwo} [rt] at 1.05 0.2
\put {$s(x)-1$} [l] at 0.03 0.2
\setsolid
\linethickness=0.9pt
\setlinear
\plot
         0         0
    0.2000         0
    0.4000    0.1100
    0.6000    0.1300
    0.8000    0.0400
    1.0000         0
/

\linethickness=0.3pt
\setdashes <3pt>
\plot
         0   0.0000
    0.0200   0.0012
    0.0400   0.0024
    0.0600   0.0039
    0.0800   0.0055
    0.1000   0.0079
    0.1200   0.0111
    0.1400   0.0148
    0.1600   0.0208
    0.1800   0.0268
    0.2000   0.0360
    0.2200   0.0456
    0.2400   0.0575
    0.2600   0.0707
    0.2800   0.0848
    0.3000   0.1007
    0.3200   0.1166
    0.3400   0.1333
    0.3600   0.1500
    0.3800   0.1657
    0.4000   0.1812
    0.4200   0.1950
    0.4400   0.2075
    0.4600   0.2188
    0.4800   0.2271
    0.5000   0.2353
    0.5200   0.2389
    0.5400   0.2424
    0.5600   0.2420
    0.5800   0.2402
    0.6000   0.2355
    0.6200   0.2275
    0.6400   0.2185
    0.6600   0.2048
    0.6800   0.1912
    0.7000   0.1746
    0.7200   0.1577
    0.7400   0.1405
    0.7600   0.1230
    0.7800   0.1062
    0.8000   0.0905
    0.8200   0.0750
    0.8400   0.0621
    0.8600   0.0493
    0.8800   0.0389
    0.9000   0.0290
    0.9200   0.0207
    0.9400   0.0138
    0.9600   0.0078
    0.9800   0.0039
    1.0000   0.0000
/

\setplotsymbol ({{\tiny $\ast$}}) \plotsymbolspacing=8pt
\setdots
\plot
         0    0.0000
    0.0200    0.0154
    0.0400    0.0308
    0.0600    0.0452
    0.0800    0.0591
    0.1000    0.0717
    0.1200    0.0827
    0.1400    0.0930
    0.1600    0.1004
    0.1800    0.1079
    0.2000    0.1125
    0.2200    0.1168
    0.2400    0.1199
    0.2600    0.1222
    0.2800    0.1243
    0.3000    0.1259
    0.3200    0.1276
    0.3400    0.1291
    0.3600    0.1306
    0.3800    0.1321
    0.4000    0.1336
    0.4200    0.1349
    0.4400    0.1359
    0.4600    0.1367
    0.4800    0.1365
    0.5000    0.1363
    0.5200    0.1342
    0.5400    0.1320
    0.5600    0.1276
    0.5800    0.1222
    0.6000    0.1148
    0.6200    0.1050
    0.6400    0.0945
    0.6600    0.0809
    0.6800    0.0673
    0.7000    0.0532
    0.7200    0.0390
    0.7400    0.0265
    0.7600    0.0150
    0.7800    0.0051
    0.8000    0.0000
    0.8200    0.0000
    0.8400    0.0000
    0.8600    0.0000
    0.8800    0.0000
    0.9000    0.0000
    0.9200    0.0000
    0.9400    0.0000
    0.9600    0.0000
    0.9800    0.0000
    1.0000    0.0000
/

\setplotsymbol ({{\tiny $\circ$}}) \plotsymbolspacing=8pt
\setdots
\plot
         0   0.0000
    0.0200   0.0036
    0.0400   0.0071
    0.0600   0.0107
    0.0800   0.0142
    0.1000   0.0179
    0.1200   0.0218
    0.1400   0.0257
    0.1600   0.0299
    0.1800   0.0340
    0.2000   0.0376
    0.2200   0.0412
    0.2400   0.0437
    0.2600   0.0496
    0.2800   0.0580
    0.3000   0.0655
    0.3200   0.0728
    0.3400   0.0778
    0.3600   0.0828
    0.3800   0.0865
    0.4000   0.0900
    0.4200   0.0932
    0.4400   0.0962
    0.4600   0.0992
    0.4800   0.1021
    0.5000   0.1050
    0.5200   0.1075
    0.5400   0.1099
    0.5600   0.1108
    0.5800   0.1110
    0.6000   0.1092
    0.6200   0.1051
    0.6400   0.1001
    0.6600   0.0914
    0.6800   0.0827
    0.7000   0.0729
    0.7200   0.0630
    0.7400   0.0535
    0.7600   0.0443
    0.7800   0.0356
    0.8000   0.0296
    0.8200   0.0286
    0.8400   0.0256
    0.8600   0.0226
    0.8800   0.0191
    0.9000   0.0155
    0.9200   0.0120
    0.9400   0.0086
    0.9600   0.0054
    0.9800   0.0027
    1.0000   0.0000
/

\setplotsymbol ({{\tiny $\bullet$}}) \plotsymbolspacing=8pt
\setdots
\plot
         0    0.0000
    0.0200    0.0000
    0.0400    0.0000
    0.0600    0.0000
    0.0800    0.0000
    0.1000    0.0000
    0.1200    0.0000
    0.1400    0.0003
    0.1600    0.0051
    0.1800    0.0099
    0.2000    0.0163
    0.2200    0.0228
    0.2400    0.0289
    0.2600    0.0388
    0.2800    0.0510
    0.3000    0.0622
    0.3200    0.0730
    0.3400    0.0804
    0.3600    0.0879
    0.3800    0.0934
    0.4000    0.0984
    0.4200    0.1031
    0.4400    0.1075
    0.4600    0.1117
    0.4800    0.1156
    0.5000    0.1194
    0.5200    0.1227
    0.5400    0.1259
    0.5600    0.1271
    0.5800    0.1276
    0.6000    0.1258
    0.6200    0.1213
    0.6400    0.1158
    0.6600    0.1061
    0.6800    0.0965
    0.7000    0.0863
    0.7200    0.0760
    0.7400    0.0662
    0.7600    0.0567
    0.7800    0.0475
    0.8000    0.0407
    0.8200    0.0388
    0.8400    0.0345
    0.8600    0.0302
    0.8800    0.0254
    0.9000    0.0205
    0.9200    0.0158
    0.9400    0.0113
    0.9600    0.0070
    0.9800    0.0035
    1.0000    0.0000
/

\endpicture
} 
\setbox\figurefive=\vbox{\hsize=\xfiglen  
\beginpicture
  \setcoordinatesystem units <\xfiglen,\yfiglen>  point at 0.0 0.0
  \setplotarea x from 0 to 1, y from 0 to 0.2
\scriptsize
  \axis bottom shiftedto y=0.0 ticks short numbered from 0 to 1 by 0.25 /
  \axis left ticks short numbered from 0.0 to 0.2 by 0.1 /
\footnotesize
\put {$\eta(x)$} [l] at 0.03 0.2
\linethickness=0.9pt
\setlinear
\setsolid
\relax
\plot
    0.00    0
    0.05    0
    0.10    0
    0.15    0
    0.20    0
    0.25    0
    0.30    0
    0.40    0.030
    0.50    0.045
    0.60    0.075
    0.70    0.150
    0.80    0.150
    0.90    0.100
    0.95    0.030
    1.00    0 /
\linethickness=0.3pt
\setdashes <3pt>
\setquadratic
\plot
         0         0
    0.0200         0
    0.0400         0
    0.0600         0
    0.0800         0
    0.1000         0
    0.1200         0
    0.1400         0
    0.1600         0
    0.1800         0
    0.2000         0
    0.2200         0
    0.2400         0
    0.2600         0
    0.2800         0
    0.3000         0
    0.3200         0
    0.3400         0
    0.3600         0
    0.3800         0
    0.4000         0
    0.4200         0
    0.4400         0
    0.4600         0
    0.4800         0
    0.5000         0
    0.5200         0
    0.5400    0.0078
    0.5600    0.0248
    0.5800    0.0422
    0.6000    0.0593
    0.6200    0.0762
    0.6400    0.0927
    0.6600    0.1075
    0.6800    0.1222
    0.7000    0.1336
    0.7200    0.1446
    0.7400    0.1523
    0.7600    0.1582
    0.7800    0.1618
    0.8000    0.1615
    0.8200    0.1605
    0.8400    0.1532
    0.8600    0.1460
    0.8800    0.1329
    0.9000    0.1186
    0.9200    0.1005
    0.9400    0.0792
    0.9600    0.0561
    0.9800    0.0281
    1.0000         0
/

\setplotsymbol ({{\tiny $\ast$}}) \plotsymbolspacing=8pt
\setdots
\plot
         0         0
    0.0200    0.0006
    0.0400    0.0013
    0.0600    0.0019
    0.0800    0.0026
    0.1000    0.0032
    0.1200    0.0039
    0.1400    0.0046
    0.1600    0.0054
    0.1800    0.0062
    0.2000    0.0071
    0.2200    0.0081
    0.2400    0.0091
    0.2600    0.0103
    0.2800    0.0116
    0.3000    0.0130
    0.3200    0.0145
    0.3400    0.0161
    0.3600    0.0177
    0.3800    0.0194
    0.4000    0.0210
    0.4200    0.0226
    0.4400    0.0241
    0.4600    0.0254
    0.4800    0.0264
    0.5000    0.0274
    0.5200    0.0276
    0.5400    0.0355
    0.5600    0.0521
    0.5800    0.0687
    0.6000    0.0846
    0.6200    0.0997
    0.6400    0.1142
    0.6600    0.1262
    0.6800    0.1383
    0.7000    0.1464
    0.7200    0.1541
    0.7400    0.1584
    0.7600    0.1608
    0.7800    0.1611
    0.8000    0.1576
    0.8200    0.1536
    0.8400    0.1443
    0.8600    0.1350
    0.8800    0.1212
    0.9000    0.1064
    0.9200    0.0889
    0.9400    0.0692
    0.9600    0.0482
    0.9800    0.0241
    1.0000         0
/

\setplotsymbol ({{\tiny $\circ$}}) \plotsymbolspacing=8pt
\setdots
\plot
         0         0
    0.0200    0.0027
    0.0400    0.0054
    0.0600    0.0081
    0.0800    0.0108
    0.1000    0.0134
    0.1200    0.0161
    0.1400    0.0187
    0.1600    0.0211
    0.1800    0.0236
    0.2000    0.0260
    0.2200    0.0283
    0.2400    0.0305
    0.2600    0.0327
    0.2800    0.0348
    0.3000    0.0367
    0.3200    0.0386
    0.3400    0.0403
    0.3600    0.0419
    0.3800    0.0432
    0.4000    0.0444
    0.4200    0.0452
    0.4400    0.0457
    0.4600    0.0459
    0.4800    0.0454
    0.5000    0.0449
    0.5200    0.0432
    0.5400    0.0492
    0.5600    0.0637
    0.5800    0.0780
    0.6000    0.0916
    0.6200    0.1042
    0.6400    0.1163
    0.6600    0.1260
    0.6800    0.1356
    0.7000    0.1416
    0.7200    0.1472
    0.7400    0.1496
    0.7600    0.1504
    0.7800    0.1492
    0.8000    0.1448
    0.8200    0.1399
    0.8400    0.1305
    0.8600    0.1211
    0.8800    0.1080
    0.9000    0.0940
    0.9200    0.0781
    0.9400    0.0604
    0.9600    0.0418
    0.9800    0.0209
    1.0000         0
/

\setplotsymbol ({{\tiny $\bullet$}}) \plotsymbolspacing=8pt
\setdots
\plot
         0         0
    0.0200    0.0028
    0.0400    0.0055
    0.0600    0.0083
    0.0800    0.0110
    0.1000    0.0137
    0.1200    0.0163
    0.1400    0.0189
    0.1600    0.0213
    0.1800    0.0237
    0.2000    0.0259
    0.2200    0.0281
    0.2400    0.0301
    0.2600    0.0320
    0.2800    0.0338
    0.3000    0.0354
    0.3200    0.0370
    0.3400    0.0383
    0.3600    0.0395
    0.3800    0.0405
    0.4000    0.0413
    0.4200    0.0418
    0.4400    0.0419
    0.4600    0.0418
    0.4800    0.0411
    0.5000    0.0403
    0.5200    0.0384
    0.5400    0.0443
    0.5600    0.0586
    0.5800    0.0728
    0.6000    0.0863
    0.6200    0.0990
    0.6400    0.1111
    0.6600    0.1208
    0.6800    0.1306
    0.7000    0.1367
    0.7200    0.1425
    0.7400    0.1453
    0.7600    0.1463
    0.7800    0.1454
    0.8000    0.1413
    0.8200    0.1368
    0.8400    0.1277
    0.8600    0.1186
    0.8800    0.1058
    0.9000    0.0922
    0.9200    0.0766
    0.9400    0.0592
    0.9600    0.0410
    0.9800    0.0205
    1.0000         0
/

\endpicture
} 
\setbox\figuresix=\vbox{\hsize=\xfiglen  
\beginpicture
  \setcoordinatesystem units <\xfiglen,\yfiglen>  point at 0.0 0.0
  \setplotarea x from 0 to 1, y from 0 to 0.2
\scriptsize
  \axis bottom shiftedto y=0.0 ticks short numbered from 0 to 1 by 0.25 /
  \axis left ticks short numbered from 0.0 to 0.2 by 0.1 /
\footnotesize
\put {\copy\figurelegendthree} [rt] at 1.05 0.2
\put {$s(x)-1$} [l] at 0.03 0.2
\setsolid
\linethickness=0.9pt
\setlinear
 \plot
         0         0
    0.2000         0
    0.4000    0.1100
    0.6000    0.1300
    0.8000    0.0400
    1.0000         0
/
\linethickness=0.3pt
\setdashes <3pt>
\plot
         0    0.0000
    0.0200    0.0000
    0.0400    0.0000
    0.0600    0.0000
    0.0800    0.0000
    0.1000    0.0000
    0.1200    0.0000
    0.1400    0.0000
    0.1600    0.0000
    0.1800    0.0000
    0.2000    0.0000
    0.2200    0.0000
    0.2400    0.0000
    0.2600    0.0000
    0.2800    0.0018
    0.3000    0.0120
    0.3200    0.0223
    0.3400    0.0342
    0.3600    0.0460
    0.3800    0.0580
    0.4000    0.0701
    0.4200    0.0813
    0.4400    0.0920
    0.4600    0.1019
    0.4800    0.1100
    0.5000    0.1180
    0.5200    0.1226
    0.5400    0.1272
    0.5600    0.1291
    0.5800    0.1300
    0.6000    0.1292
    0.6200    0.1263
    0.6400    0.1229
    0.6600    0.1171
    0.6800    0.1112
    0.7000    0.1037
    0.7200    0.0959
    0.7400    0.0878
    0.7600    0.0793
    0.7800    0.0710
    0.8000    0.0628
    0.8200    0.0547
    0.8400    0.0472
    0.8600    0.0398
    0.8800    0.0332
    0.9000    0.0267
    0.9200    0.0207
    0.9400    0.0150
    0.9600    0.0096
    0.9800    0.0048
    1.0000    0.0000
/

\setplotsymbol ({{\tiny $\ast$}}) \plotsymbolspacing=8pt
\setdots
\plot
         0    0.0000
    0.0200    0.0049
    0.0400    0.0098
    0.0600    0.0145
    0.0800    0.0192
    0.1000    0.0238
    0.1200    0.0282
    0.1400    0.0324
    0.1600    0.0362
    0.1800    0.0400
    0.2000    0.0430
    0.2200    0.0460
    0.2400    0.0483
    0.2600    0.0501
    0.2800    0.0533
    0.3000    0.0640
    0.3200    0.0747
    0.3400    0.0855
    0.3600    0.0964
    0.3800    0.1062
    0.4000    0.1157
    0.4200    0.1240
    0.4400    0.1311
    0.4600    0.1373
    0.4800    0.1413
    0.5000    0.1452
    0.5200    0.1458
    0.5400    0.1464
    0.5600    0.1445
    0.5800    0.1417
    0.6000    0.1374
    0.6200    0.1313
    0.6400    0.1248
    0.6600    0.1164
    0.6800    0.1079
    0.7000    0.0985
    0.7200    0.0890
    0.7400    0.0797
    0.7600    0.0705
    0.7800    0.0618
    0.8000    0.0538
    0.8200    0.0459
    0.8400    0.0394
    0.8600    0.0329
    0.8800    0.0274
    0.9000    0.0221
    0.9200    0.0172
    0.9400    0.0126
    0.9600    0.0081
    0.9800    0.0041
    1.0000    0.0000
/

\setplotsymbol ({{\tiny $\circ$}}) \plotsymbolspacing=8pt
\setdots
\plot
         0    0.0000
    0.0200    0.0055
    0.0400    0.0110
    0.0600    0.0163
    0.0800    0.0215
    0.1000    0.0264
    0.1200    0.0311
    0.1400    0.0355
    0.1600    0.0392
    0.1800    0.0430
    0.2000    0.0457
    0.2200    0.0483
    0.2400    0.0501
    0.2600    0.0513
    0.2800    0.0538
    0.3000    0.0636
    0.3200    0.0735
    0.3400    0.0832
    0.3600    0.0930
    0.3800    0.1016
    0.4000    0.1100
    0.4200    0.1170
    0.4400    0.1230
    0.4600    0.1282
    0.4800    0.1313
    0.5000    0.1344
    0.5200    0.1346
    0.5400    0.1348
    0.5600    0.1328
    0.5800    0.1300
    0.6000    0.1258
    0.6200    0.1199
    0.6400    0.1136
    0.6600    0.1054
    0.6800    0.0972
    0.7000    0.0882
    0.7200    0.0792
    0.7400    0.0705
    0.7600    0.0621
    0.7800    0.0542
    0.8000    0.0472
    0.8200    0.0404
    0.8400    0.0350
    0.8600    0.0297
    0.8800    0.0251
    0.9000    0.0208
    0.9200    0.0166
    0.9400    0.0124
    0.9600    0.0083
    0.9800    0.0042
    1.0000    0.0000
/

\setplotsymbol ({{\tiny $\bullet$}}) \plotsymbolspacing=8pt
\setdots
\plot
         0    0.0000
    0.0200    0.0047
    0.0400    0.0093
    0.0600    0.0138
    0.0800    0.0182
    0.1000    0.0222
    0.1200    0.0260
    0.1400    0.0296
    0.1600    0.0326
    0.1800    0.0355
    0.2000    0.0377
    0.2200    0.0398
    0.2400    0.0411
    0.2600    0.0421
    0.2800    0.0444
    0.3000    0.0542
    0.3200    0.0640
    0.3400    0.0739
    0.3600    0.0837
    0.3800    0.0925
    0.4000    0.1010
    0.4200    0.1083
    0.4400    0.1145
    0.4600    0.1200
    0.4800    0.1234
    0.5000    0.1269
    0.5200    0.1277
    0.5400    0.1285
    0.5600    0.1271
    0.5800    0.1250
    0.6000    0.1213
    0.6200    0.1160
    0.6400    0.1102
    0.6600    0.1022
    0.6800    0.0943
    0.7000    0.0855
    0.7200    0.0766
    0.7400    0.0682
    0.7600    0.0600
    0.7800    0.0523
    0.8000    0.0457
    0.8200    0.0393
    0.8400    0.0343
    0.8600    0.0294
    0.8800    0.0252
    0.9000    0.0212
    0.9200    0.0171
    0.9400    0.0130
    0.9600    0.0089
    0.9800    0.0044
    1.0000    0.0000
/

\endpicture
} 
\setbox\figureseven=\vbox{\hsize=\xfiglen  
\beginpicture
  \setcoordinatesystem units <\xfiglen,\yfiglen>  point at 0.0 0.0
  \setplotarea x from 0 to 1, y from 0 to 0.2
\scriptsize
  \axis bottom shiftedto y=0.0 ticks short numbered from 0 to 1 by 0.25 /
  \axis left ticks short numbered from 0.0 to 0.2 by 0.1 /
\footnotesize
\put {$\eta(x)$} [l] at 0.03 0.2
\linethickness=0.9pt
\setlinear
\setsolid
\relax
\plot
         0         0
    0.0500         0
    0.1000    0.0500
    0.1500    0.1000
    0.2000    0.1300
    0.2500    0.1300
    0.3000    0.0800
    0.4000    0.0300
    0.5000    0.0450
    0.6000    0.0750
    0.7000    0.1500
    0.8000    0.1500
    0.9000    0.1000
    0.9500         0
    1.0000         0
/

\linethickness=0.3pt
\setdashes <3pt>
\setquadratic
\plot
         0         0
    0.0200    0.0001
    0.0400    0.0001
    0.0600    0.0003
    0.0800    0.0004
    0.1000    0.0008
    0.1200    0.0014
    0.1400    0.0022
    0.1600    0.0035
    0.1800    0.0048
    0.2000    0.0063
    0.2200    0.0079
    0.2400    0.0092
    0.2600    0.0104
    0.2800    0.0112
    0.3000    0.0114
    0.3200    0.0116
    0.3400    0.0107
    0.3600    0.0099
    0.3800    0.0087
    0.4000    0.0074
    0.4200    0.0066
    0.4400    0.0061
    0.4600    0.0061
    0.4800    0.0076
    0.5000    0.0091
    0.5200    0.0137
    0.5400    0.0183
    0.5600    0.0253
    0.5800    0.0332
    0.6000    0.0424
    0.6200    0.0532
    0.6400    0.0643
    0.6600    0.0766
    0.6800    0.0888
    0.7000    0.1010
    0.7200    0.1131
    0.7400    0.1242
    0.7600    0.1347
    0.7800    0.1440
    0.8000    0.1511
    0.8200    0.1577
    0.8400    0.1591
    0.8600    0.1604
    0.8800    0.1541
    0.9000    0.1461
    0.9200    0.1304
    0.9400    0.1082
    0.9600    0.0811
    0.9800    0.0405
    1.0000         0
/

\setplotsymbol ({{\tiny $\ast$}}) \plotsymbolspacing=8pt
\setdots
\plot
         0         0
    0.0200    0.0116
    0.0400    0.0233
    0.0600    0.0348
    0.0800    0.0463
    0.1000    0.0573
    0.1200    0.0679
    0.1400    0.0780
    0.1600    0.0865
    0.1800    0.0949
    0.2000    0.1001
    0.2200    0.1050
    0.2400    0.1069
    0.2600    0.1071
    0.2800    0.1052
    0.3000    0.1000
    0.3200    0.0943
    0.3400    0.0840
    0.3600    0.0738
    0.3800    0.0707
    0.4000    0.0675
    0.4200    0.0634
    0.4400    0.0585
    0.4600    0.0540
    0.4800    0.0504
    0.5000    0.0468
    0.5200    0.0473
    0.5400    0.0479
    0.5600    0.0525
    0.5800    0.0585
    0.6000    0.0666
    0.6200    0.0771
    0.6400    0.0878
    0.6600    0.0997
    0.6800    0.1115
    0.7000    0.1213
    0.7200    0.1308
    0.7400    0.1370
    0.7600    0.1412
    0.7800    0.1431
    0.8000    0.1406
    0.8200    0.1375
    0.8400    0.1284
    0.8600    0.1194
    0.8800    0.1061
    0.9000    0.0919
    0.9200    0.0758
    0.9400    0.0582
    0.9600    0.0399
    0.9800    0.0200
    1.0000         0
/

\setplotsymbol ({{\tiny $\circ$}}) \plotsymbolspacing=8pt
\setdots
\plot
         0         0
    0.0200    0.0022
    0.0400    0.0044
    0.0600    0.0066
    0.0800    0.0089
    0.1000    0.0115
    0.1200    0.0144
    0.1400    0.0174
    0.1600    0.0213
    0.1800    0.0251
    0.2000    0.0295
    0.2200    0.0339
    0.2400    0.0378
    0.2600    0.0415
    0.2800    0.0442
    0.3000    0.0453
    0.3200    0.0461
    0.3400    0.0434
    0.3600    0.0407
    0.3800    0.0448
    0.4000    0.0488
    0.4200    0.0512
    0.4400    0.0522
    0.4600    0.0532
    0.4800    0.0541
    0.5000    0.0551
    0.5200    0.0590
    0.5400    0.0630
    0.5600    0.0702
    0.5800    0.0787
    0.6000    0.0886
    0.6200    0.1003
    0.6400    0.1121
    0.6600    0.1237
    0.6800    0.1353
    0.7000    0.1432
    0.7200    0.1508
    0.7400    0.1542
    0.7600    0.1551
    0.7800    0.1532
    0.8000    0.1464
    0.8200    0.1391
    0.8400    0.1258
    0.8600    0.1125
    0.8800    0.0964
    0.9000    0.0796
    0.9200    0.0626
    0.9400    0.0457
    0.9600    0.0293
    0.9800    0.0147
    1.0000         0
/

\setplotsymbol ({{\tiny $\bullet$}}) \plotsymbolspacing=8pt
\setdots
\plot
         0         0
    0.0200         0
    0.0400         0
    0.0600         0
    0.0800         0
    0.1000         0
    0.1200         0
    0.1400         0
    0.1600    0.0020
    0.1800    0.0049
    0.2000    0.0095
    0.2200    0.0143
    0.2400    0.0192
    0.2600    0.0242
    0.2800    0.0284
    0.3000    0.0312
    0.3200    0.0336
    0.3400    0.0320
    0.3600    0.0304
    0.3800    0.0349
    0.4000    0.0391
    0.4200    0.0413
    0.4400    0.0419
    0.4600    0.0425
    0.4800    0.0431
    0.5000    0.0436
    0.5200    0.0478
    0.5400    0.0519
    0.5600    0.0599
    0.5800    0.0694
    0.6000    0.0805
    0.6200    0.0939
    0.6400    0.1072
    0.6600    0.1205
    0.6800    0.1338
    0.7000    0.1434
    0.7200    0.1526
    0.7400    0.1573
    0.7600    0.1595
    0.7800    0.1588
    0.8000    0.1530
    0.8200    0.1465
    0.8400    0.1337
    0.8600    0.1209
    0.8800    0.1046
    0.9000    0.0876
    0.9200    0.0700
    0.9400    0.0518
    0.9600    0.0340
    0.9800    0.0170
    1.0000         0
/

\endpicture
} 
\setbox\figureeight=\vbox{\hsize=\xfiglen  
\beginpicture
  \setcoordinatesystem units <\xfiglen,\yfiglen>  point at 0.0 0.0
  \setplotarea x from 0 to 1, y from 0 to 0.2
\scriptsize
  \axis bottom shiftedto y=0.0 ticks short numbered from 0 to 1 by 0.25 /
  \axis left ticks short numbered from 0.0 to 0.2 by 0.1 /
\footnotesize
\put {\copy\figurelegendfour} [rt] at 1.05 0.2
\put {$s(x)-1$} [l] at 0.03 0.2
\setsolid
\linethickness=0.9pt
\setlinear
 \plot
         0         0
    0.2000         0
    0.4000    0.1100
    0.6000    0.1300
    0.8000    0.0400
    1.0000         0
/
\linethickness=0.3pt
\setdashes <3pt>
\plot
         0    0.0000
    0.0200    0.0178
    0.0400    0.0355
    0.0600    0.0520
    0.0800    0.0680
    0.1000    0.0822
    0.1200    0.0943
    0.1400    0.1056
    0.1600    0.1133
    0.1800    0.1209
    0.2000    0.1250
    0.2200    0.1287
    0.2400    0.1311
    0.2600    0.1326
    0.2800    0.1340
    0.3000    0.1351
    0.3200    0.1363
    0.3400    0.1379
    0.3600    0.1395
    0.3800    0.1415
    0.4000    0.1436
    0.4200    0.1455
    0.4400    0.1472
    0.4600    0.1485
    0.4800    0.1486
    0.5000    0.1488
    0.5200    0.1463
    0.5400    0.1439
    0.5600    0.1387
    0.5800    0.1325
    0.6000    0.1242
    0.6200    0.1134
    0.6400    0.1018
    0.6600    0.0871
    0.6800    0.0725
    0.7000    0.0575
    0.7200    0.0425
    0.7400    0.0294
    0.7600    0.0174
    0.7800    0.0073
    0.8000    0.0005
    0.8200    0.0000
    0.8400    0.0000
    0.8600    0.0000
    0.8800    0.0000
    0.9000    0.0000
    0.9200    0.0000
    0.9400    0.0000
    0.9600    0.0000
    0.9800    0.0000
    1.0000    0.0000
/

\setplotsymbol ({{\tiny $\ast$}}) \plotsymbolspacing=8pt
\setdots
\plot
         0    0.0000
    0.0200    0.0072
    0.0400    0.0144
    0.0600    0.0214
    0.0800    0.0283
    0.1000    0.0349
    0.1200    0.0411
    0.1400    0.0472
    0.1600    0.0526
    0.1800    0.0579
    0.2000    0.0623
    0.2200    0.0666
    0.2400    0.0700
    0.2600    0.0730
    0.2800    0.0756
    0.3000    0.0824
    0.3200    0.0925
    0.3400    0.1017
    0.3600    0.1110
    0.3800    0.1194
    0.4000    0.1275
    0.4200    0.1349
    0.4400    0.1416
    0.4600    0.1476
    0.4800    0.1521
    0.5000    0.1566
    0.5200    0.1583
    0.5400    0.1601
    0.5600    0.1590
    0.5800    0.1569
    0.6000    0.1524
    0.6200    0.1451
    0.6400    0.1370
    0.6600    0.1251
    0.6800    0.1132
    0.7000    0.1002
    0.7200    0.0871
    0.7400    0.0749
    0.7600    0.0631
    0.7800    0.0523
    0.8000    0.0431
    0.8200    0.0401
    0.8400    0.0358
    0.8600    0.0315
    0.8800    0.0266
    0.9000    0.0216
    0.9200    0.0167
    0.9400    0.0119
    0.9600    0.0074
    0.9800    0.0037
    1.0000    0.0000
/

\setplotsymbol ({{\tiny $\circ$}}) \plotsymbolspacing=8pt
\setdots
\plot
         0    0.0000
    0.0200    0.0017
    0.0400    0.0034
    0.0600    0.0053
    0.0800    0.0073
    0.1000    0.0099
    0.1200    0.0131
    0.1400    0.0167
    0.1600    0.0217
    0.1800    0.0267
    0.2000    0.0323
    0.2200    0.0379
    0.2400    0.0430
    0.2600    0.0477
    0.2800    0.0516
    0.3000    0.0592
    0.3200    0.0700
    0.3400    0.0787
    0.3600    0.0874
    0.3800    0.0948
    0.4000    0.1018
    0.4200    0.1082
    0.4400    0.1141
    0.4600    0.1197
    0.4800    0.1244
    0.5000    0.1291
    0.5200    0.1324
    0.5400    0.1356
    0.5600    0.1368
    0.5800    0.1371
    0.6000    0.1354
    0.6200    0.1311
    0.6400    0.1259
    0.6600    0.1169
    0.6800    0.1079
    0.7000    0.0979
    0.7200    0.0879
    0.7400    0.0781
    0.7600    0.0686
    0.7800    0.0595
    0.8000    0.0513
    0.8200    0.0490
    0.8400    0.0445
    0.8600    0.0399
    0.8800    0.0342
    0.9000    0.0282
    0.9200    0.0222
    0.9400    0.0160
    0.9600    0.0101
    0.9800    0.0051
    1.0000    0.0000
/

\setplotsymbol ({{\tiny $\bullet$}}) \plotsymbolspacing=8pt
\setdots
\plot
         0    0.0000
    0.0200    0.0010
    0.0400    0.0020
    0.0600    0.0032
    0.0800    0.0045
    0.1000    0.0065
    0.1200    0.0095
    0.1400    0.0130
    0.1600    0.0186
    0.1800    0.0243
    0.2000    0.0310
    0.2200    0.0378
    0.2400    0.0439
    0.2600    0.0497
    0.2800    0.0545
    0.3000    0.0625
    0.3200    0.0736
    0.3400    0.0819
    0.3600    0.0902
    0.3800    0.0967
    0.4000    0.1028
    0.4200    0.1082
    0.4400    0.1131
    0.4600    0.1176
    0.4800    0.1213
    0.5000    0.1249
    0.5200    0.1273
    0.5400    0.1297
    0.5600    0.1302
    0.5800    0.1299
    0.6000    0.1275
    0.6200    0.1226
    0.6400    0.1169
    0.6600    0.1074
    0.6800    0.0979
    0.7000    0.0881
    0.7200    0.0782
    0.7400    0.0690
    0.7600    0.0601
    0.7800    0.0516
    0.8000    0.0439
    0.8200    0.0423
    0.8400    0.0385
    0.8600    0.0346
    0.8800    0.0298
    0.9000    0.0247
    0.9200    0.0194
    0.9400    0.0140
    0.9600    0.0088
    0.9800    0.0044
    1.0000    0.0000
/

\endpicture
} 
\setbox\figurenine=\vbox{
\beginpicture
  \setcoordinatesystem units <\xfiglen,\yfiglen>  point at 0.0 0.0
  \setplotarea x from 0 to 30, y from 0.5 to -5
\scriptsize
  \axis bottom shiftedto y=-5.0 ticks short numbered from 0 to 30 by 5 /
  \axis left ticks short numbered from -5 to 0 by 1 /
\footnotesize
\put {$\log_{10}\sigma_n$} [lt] at 0.7 0.7
\put {$n$} [rb] at 30 -4.85
 \multiput {$\bullet$} at
    1.0000    0.0787
    2.0000   -0.7060
    3.0000   -1.1187
    4.0000   -1.3377
    5.0000   -1.4964
    6.0000   -1.6221
    7.0000   -1.7305
    8.0000   -1.8277
    9.0000   -1.9182
   10.0000   -2.0038
   11.0000   -2.0863
   12.0000   -2.1664
   13.0000   -2.2447
   14.0000   -2.3221
   15.0000   -2.4014
   16.0000   -2.4876
   17.0000   -2.5853
   18.0000   -2.7008
   19.0000   -2.8420
   20.0000   -3.0206
   21.0000   -3.2598
   22.0000   -3.3587
   23.0000   -3.5853
   24.0000   -3.6378
   25.0000   -3.7524
   26.0000   -3.9024
   27.0000   -4.0444
   28.0000   -4.1831
   29.0000   -4.3208
   30.0000   -4.4593 /

\setsolid
\linethickness=0.9pt
\setlinear
\endpicture
} 

\begin{document}
\title{
\revision{
On the simultaneous reconstruction of the nonlinearity coefficient and the sound speed in the Westervelt equation
}
}
\author{Barbara Kaltenbacher\footnote{
Department of Mathematics,
Alpen-Adria-Universit\"at Klagenfurt.
barbara.kaltenbacher@aau.at.}
\and
William Rundell\footnote{
Department of Mathematics,
Texas A\&M University,
Texas 77843. 
rundell@math.tamu.edu}
}
\date{\vskip-3ex}
\maketitle  
\vspace{-10pt}
\begin{abstract}
This paper considers the Westervelt equation, one of the most widely used
models in nonlinear acoustics, and seeks to recover two
spatially-dependent parameters of physical importance from time-trace boundary measurements.
Specifically, these are the nonlinearity parameter $\kappa(x)$ often
referred to as $B/A$ in the acoustics literature and the wave speed $c_0(x)$.
The determination of the spatial change in these quantities can be used as a means of imaging. 
We consider identifiability from one or two boundary measurements as relevant in these applications. 
\anapap{
More precisely, we provide results on local uniqueness of $\kappa(x)$ from a single observation and on simultaneous identifiability of $\kappa(x)$ and $c_0(x)$ from two measurements.
}   
For a reformulation of the problem in terms of the squared slowness $\ssl=1/c_0^2$ and the combined coefficient $\nlc=\frac{\kappa}{c_0^2}$ we devise a frozen Newton method and prove its convergence.
The effectiveness (and limitations) of this iterative scheme are demonstrated by numerical examples.
\end{abstract}

\leftline{\small \qquad\qquad
}
\smallskip
\leftline{\small \qquad\qquad 
\textbf{{\textsc AMS}} {\bf classification:} 35R30, 35R11, 35L70, 78A46}

\section{Introduction}
Imaging with ultrasound has a long and successful history 
based on a vast range of applications. 
However, as is often the case, the use of lower frequencies naturally leads to
lower resolution and at higher frequencies sound propagation is affected by
scattering and stronger attenuation.
Enhanced ultrasound-based techniques such as 
nonlinearity parameter imaging
\cite{Bjorno1986, BurovGurinovichRudenkoTagunov1994, Cain1986, 
IchidaSatoLinzer1983, PanfilovaSlounWijkstraSapozhnikovMischi:2021, VarrayBassetTortoliCachard2011,ZhangChenGong2001, ZhangChenYe1996}, 
harmonic imaging 
\cite{AnvariForsbergSamir15,Uppal2010,VarrayBassetTortoliCachard2011}, 
and vibro-acoustography 
\cite{FatemiGreenleaf1998,FatemiGreenleaf1999,vibroacoustics,Malcolmetal2007,Malcolmetal2008}
have been developed to overcome these drawbacks and improve imaging quality.
They make use of {\it nonlinear effects\/} that arise at higher intensities
or when waves interact and are characterised by a multiplicative coefficient
that is usually called parameter of nonlinearity and denoted by $B/A$.
We will here use the mathematically convenient
abbreviation $\kappaintro$
for a quantity containing $B/A$.
This coefficient depends on tissue properties and therefore varies in the
spatial direction, $\kappaintro=\kappaintro(x)$.

While ultrasound imaging relies on the propagation of sound waves and is
therefore physically and mathematically correctly described by some wave-type
partial differential equation ({\sc pde}),
algorithms implemented in modern ultrasound scanners make use of model
simplifications that allow one to apply methods from signal processing
(beamforming, filtering) to generate an image based on the principles of
transmission and reflection, based on differences in the acoustic impedance $Z$. 
These simplifications are not able to capture nonlinearity so that one has to
return to the {\sc pde} model and consider $\kappaintro=\kappaintro(x)$
(and often also the speed of sound $c_0=c_0(x)$) as a spatially variable
 coefficient.

\anapap{
A novel approach of ultrasound tomography
\cite{AlsakerCardenasFuruieMueller2021,Gemmeke_etal:2017,Greenleaf_etal:1974,JavaherianLuckaCox20,LuckaPerezlivaTreebyCox2022,MuellerCardenasFuruie2021}
aims at recovering the speed of sound $c_0$ as a space dependent function in view of the fact that the impedance $Z$ 
responsible for transmission and reflection is tied to $c_0$ via the
identity $Z=\varrho_0\, c_0$, where $\varrho_0$ is the mass density.
Likewise, nonlinearity parameter imaging needs to be viewed as the identification of a space-dependent coefficient $\kappaintro=\kappaintro(x)$ in a {\sc pde} and as such
it will often arise alongside together with the recovery of the spatially
varying sound speed $c_0=c_0(x)$. 
}

In the following subsections we provide more background on the mathematical
models. In particular we will  show at which position in the {\sc pde} these
coefficients appear, which of course is a factor crucial for their recovery.
We then describe the inverse problem and the basic method of its solution.

\bigskip

We consider, as one of the most established classical model of nonlinear acoustics, the Westervelt equation in pressure formulation
\begin{equation}\label{Westervelt}
u_{tt}-c_0^2\Delta u - b\Delta u_t = \kappaintro(u^2)_{tt} + r
\mbox{ in }(0,T)\times\Omega\,, 
\end{equation}
where $u$ is the acoustic pressure, $c_0$ the speed of sound, $b$ the diffusivity of sound, $\varrho_0$ the mass density, and $\kappaintro$ 
$=\frac{\beta_a}{\varrho_0 c_0^2}=\tfrac{1}{\varrho_0 c_0^2}(\tfrac{B}{2A}+1)$
contains the nonlinearity parameter $\beta_a$ or $B/A$. 
\footnote{
More precisely, the {\sc pde} is 
$\frac{1}{\lambda(x)}u_{tt}-\nabla\cdot(\frac{1}{\varrho_0(x)}\nabla u) - b D u = \kappaintro(u^2)_{tt}+r$ with $u$ being the pressure, $\lambda$ the bulk modulus, $\varrho_0$ the mass density, and $c_0=\frac{\lambda}{\varrho_0}$ the sound speed, (cf., e.g., \cite{BambergerGlowinskiTran,LuckaPerezlivaTreebyCox2022} for the linear case). Spatial variability of $\varrho$ is not relevant in our context; rather, dependence of $c_0$ on $x$ is due to variability of $\lambda$.
}
We assume \eqref{Westervelt} to hold in a domain $\Omega\subseteq\mathbb{R}^3$ and equip it with initial conditions $u(t=0)=u_0$, $u_t(t=0)=u_1$, as well as
absorbing or impedance  boundary conditions on the rest of the boundary
to enable restriction to a bounded computational domain $\Omega$, which without loss of generality we can assume to be smooth.
The space- and time-dependent interior source term $r$ in \eqref{Westervelt} models
excitation by a piezoelectric transducer array. 

The pressure data $h$ taken at the receiver array 
is expressed as a Dirichlet trace on some manifold $\Sigma$ immersed in the computational domain $\Omega$ or attached to its boundary $\Sigma\in\overline{\Omega}$
\begin{equation}\label{observation}
h(t,x) = u(t,x), \quad(t,x)\in(0,T)\times \Sigma.
\end{equation}
Note that our setting also allows $\Sigma$ to be a subset of discrete points on a manifold.

The inverse problem of nonlinearity parameter tomography consists of reconstructing $\kappaintro=\kappaintro(x)$ 
from measurements \eqref{observation}. 
Often, the speed of sound is a function of the space variables as well $c_0=c_0(x)$ and needs to be recovered alongside with $\kappa$.
This is a natural requirement as the sound speed will vary between objects to be imaged and also from the background.

Reconstruction of $c_0=c_0(x)$ as a \textsc{pde} coefficient is actually already being done in ultrasound tomography
\cite{AlsakerCardenasFuruieMueller2021,Gemmeke_etal:2017,Greenleaf_etal:1974,JavaherianLuckaCox20,MuellerCardenasFuruie2021}
but in a linear wave equation, that is, with $\kappaintro=0$.
We mention in passing that in principle the mass density $\varrho_0$ also varies
in the spatial direction.
However, in ultrasound imaging, this coefficient  does not play a significant
role and is therefore usually neglected.\relax
\footnote{The notation $c_0$, $\varrho_0$ for the (reference) sound speed and mass density, respectively, refers to he usual decomposition of the mass density into a reference and a fluctuation part $\varrho = \varrho_0+\varrho_\sim$, and correspondingly for $c_0=\frac{\lambda}{\varrho_0}$, where $\lambda$ is the bulk modulus. While the total mass density $\varrho$ would be subject to a balance law (namely conservation of mass) and thus appear as one of the states in a \textsc{PDE} model, its appearance as a coefficient only affects the reference part $\varrho_0$. This is due to the typical expansion rules (known as Blackstock's scheme in the nonlinear acoustics literature) for obtaining linear and quadratic acoustic wave equations from nonlinear balance and constitutive laws, cf. e.g., \cite{HamiltonBlackstock97,reviewNonlinearAcoustics} and the references therein.
}

We refer to \cite{AcostaUhlmannZhai:2022,
nonlinearity_imaging_Westervelt,nonlinearity_imaging_fracWest,nonlinearityimaging} for results related to the identification of the nonlinearity coefficient $\kappa$ alone.
In \cite{AcostaUhlmannZhai:2022} its uniqueness from the whole Neumann-Dirichlet map (instead of the single measurement \eqref{observation}) is shown; \cite{nonlinearityimaging} provides a uniqueness and conditional stability result for the linearised problem of identifying $\kappa$ in a higher order model of nonlinear acoustics in place of the  Westervelt equation.
In \cite{nonlinearity_imaging_Westervelt,nonlinearity_imaging_fracWest} we have proven injectivity of the linearised forward operator mapping $\kappa$ to $h$ in the Westervelt equation with classical strong damping and also with some fractional damping models as relevant in ultrasonics. 
\anapap{
This will serve as a basis for proving local uniqueness of $\kappa$ by means of the Inverse Function Theorem.
}

\anapap{
Besides this, the aim of this paper is to study simultaneous identification of $\kappa$ and $c_0$ as space variable functions. Indeed, 
we will show that $\kappa(x)$ is locally unique even for unknown $c_0(x)$ in general space dimension $d\in\{1,2,3\}$ in Section~\ref{sec:uniqueness_kappa} and provide a result on simultaneous identifiability of $\kappa(x)$ and $c_0(x)$ in one space dimension based on inverse Sturm-Liouville theory in Section~\ref{sec:uniqueness_kappa_c0_1-d}. Moreover,
} 

The aim of this paper is to provide results on the {\em simultaneous} recovery of $\kappa(x)$ and $c_0(x)$.
In Section~\ref{sec:uniqueness_kappa_c0_lin} we will prove injectivity of the linearised forward operator from measurements with two excitations.
This serves as a basis for applying a frozen Newton method and showing its convergence in Section~\ref{sec:Newton}. 
Numerical reconstruction results are provided in Section~\ref{sec:reconstructions}.

\subsection{The inverse problem}\label{sec:inverse}
Consider identification of the space dependent nonlinearity coefficient $\kappa(x)$ and sound speed $c_0(x)$ 
for the attenuated Westervelt equation in pressure form
\begin{equation}\label{eqn:Westervelt_init_D_intro}
\begin{aligned}
&\bigl(u-\kappa(x)u^2\bigr)_{tt}-c_0(x)^2\triangle u + D[u] = r \quad
\mbox{ in }\Omega\times(0,T)\\
\partial_\nu u+\gamma u&=0 \mbox{ on }\partial\Omega\times(0,T),\quad
u(0)=0, \quad u_t(0)=0 \quad \mbox{ in }\;\Omega
\end{aligned}
\end{equation}
from observations of the acoustic pressure 
\begin{equation}\label{eqn:obs}
h(x,t)=u(x,t) \,,\quad x\in\Sigma\,, \quad t\in(0,T).
\end{equation}
The physical meanings of the quantities in this model are listed in Table~\ref{tab:physics}, where 
we assume $\varrho_0>0$, $b>0$, $\gamma\geq0$ to be known constants, whereas $B/A$ (and therefore $\kappa$), as well as $c_0$ may depend on the $x$ variables.
\begin{table}[ht]
\begin{center}
\begin{tabular}{l}
$u\ \ldots\ $ pressure $\ [g m^{-1} s^{-2}]$\\
$c_0\ \ldots\ $ sound speed $\ [m\, s^{-1}]$\\
$\kappa=\frac{B/A+2}{\varrho_0 c_0^2}\ \ldots\ $ nonlinearity coefficient $\ [g^{-1}\, m \, s^2 ]$\\
$\varrho_0\in\mathbb{R}^+\ \ldots\ $ mass density $\ [g\, m^{-3}]$\\
$B/A\ \ldots\ $ nonlinearity parameter $\ [1]$\\
$b\in\mathbb{R}^+\ \ldots\ $ diffusivity of sound $\ [m^2\, s^{-1}]$\\
$\gamma\in\mathbb{R}_0^+\ \ldots\ $ boundary impedance $\ [m^{-1}]$\\
\end{tabular}
\caption{Physical quantities appearing in the {\sc pde}s. \label{tab:physics}}
\end{center}
\end{table}

In equation \eqref{eqn:Westervelt_init_D_intro}, the damping term $D$ is defined by one of the two following fractional damping models 
$$
\begin{aligned}
D &= b (-\triangle)^\beta \partial_t^\alpha 
\quad &\mbox{ (combination of Caputo-Wismer-Kelvin and Chen-Holm -- {\CWKCH})}
\\
D &= b_1 (-\triangle) \partial_t^{\alpha_1}  + b_2 \partial_t^{\alpha_2+2} 
 &\ \mbox{ (fractional Zener -- \FZ)}
\end{aligned}
$$
(for more details see, e.g., \cite{nonlinearity_imaging_fracWest} and the references therein, in particular
\cite{Caputo:1967,Wismer:2006,ChenHolm:2004} for \CWKCH and \cite{HolmNaesholm:2011,Mainardi:2010,CaiChenFangHolm_survey2018} for \FZ).

The time fractional derivatives appearing in the damping models are defined by the Djrbashian-Caputo derivative
\[
\partial_t^\alpha u =I^{1-\alpha} u_t 
\]
with the Abel integral operator 
\[
(I^{1-\alpha} v)(t)  \frac{1}{\Gamma(1-\alpha)}\int_0^t \frac{v(s)}{(t-s)^{\alpha}}\,ds
\]
and $\alpha\in(0,1)$.
For defining fractional powers of the negative Laplacian $-\triangle$ with impedance boundary conditions in the \CWKCH case, we use the spectral definition
\[
((-\triangle)^\beta v)(x) = \sum_{j=1}^\infty \lambda_j^\beta\sum_{k\in K^{\lambda_j}} \langle v,\varphi_j^k\rangle \varphi_j^k(x).
\]

Excitation is modeled by an interior space and time dependent source term $r$, which indeed allows to describe the acoustic signal emitted by a transducer array immersed in the domain $\Omega$, see also \cite{periodicWestervelt}.

\anapap{
Equation \eqref{eqn:Westervelt_init_D_intro} is the formulation we will use for the uniqueness proofs, where it is convenient to attach $c_0^2(x)$ to $-\triangle$ within $\mathcal{A}$ in order to be able to apply inverse Sturm-Liouville theory.
For these proofs, the inverse problem will make use of the forward operator $\mathbb{F}$ that takes $\kappa$ to the residues of $\widehat{\textup{tr}_\Sigma u}$ (the Laplace transform of the trace of $u$ on the observation set $\Sigma$) at the poles of the Laplace transformed observation $\hat{h}$.
In order to be able to take Laplace transforms, throughout Section~\ref{sec:uniqueness}  we will assume to have observations on the whole positive timeline, which in case of $r$ being analytic with respect to time (for example, just vanishing) from a time instance $\underline{T}$ on, follows by analytic continuation of the Fourier components of $u$, which will be shown to satisfy linear constant coefficient {\sc ode}s from a time instance $T_*$ on in Section~\ref{sec:uniqueness}.
}

\bigskip

In most of this paper, we will work with the following alternative formulation that moves the spatially variable coefficient $c_0(x)$ away from the Laplacian and thus leads to a symmetric positive (as well as relatively simple) elliptic differential operator in the equation.
\footnote{An alternative to achieve symmetry would be to use the weighted $L^2$ inner product with weight function $1/c_0^2(x)$}
To this end, we divide \eqref{eqn:Westervelt_init_D_intro} by $c_0^2(x)$ and rewrite it, using the new coefficient functions $\ssl(x)$, $\nlc(x)$, as
\begin{equation}\label{eqn:Westervelt_ssl_nlc_intro}
\begin{aligned}
&\bigl(\ssl(x) u-\nlc(x)u^2\bigr)_{tt}-\triangle u + \tilde{D} u = \tilde{r} \quad
&\mbox{ in }\Omega\times(0,T)\\
&\partial_\nu u+\gamma u=0 \mbox{ on }\partial\Omega\times(0,T),\quad
u(0)=0, \quad u_t(0)=0 \quad &\mbox{ in }\;\Omega
\end{aligned}
\end{equation}
where (again with physical units in brackets)\hfill\break
\phantom{grrr} $\ssl=\frac{1}{c_0^2} \ldots\ $
squared slowness $\ [m^{-2} s^2]$\hfill\break
\phantom{grrr}
$\nlc=\frac{\kappa}{c_0^2}=\frac{\beta_a}{\varrho_0 c_0^4}=\frac{B/A+2}{\varrho_0 c_0^4}=\frac{B/A+2}{\varrho_0}\ssl^2\ \ldots\ $
nonlinearity coefficient $\ [g^{-1} m^{-1} s^4]$\hfill\break
\phantom{grrr}
$\tilde{D}=\tilde{b} (-\triangle)^\beta \partial_t^\alpha \mbox{ ({\CWKCH}) }
\quad \mbox{ or }\quad
\tilde{D}=\tilde{b}_1 (-\triangle) \partial_t^{\alpha_1}  + \tilde{b}_2 \partial_t^{\alpha_2+2} 
\mbox{ (\FZ)}.
$ \hfill\break
%
Note that we neglect variability of the damping and driving terms term with division by $c_0(x)^2$ and assume $\tilde{D}$ to come with constant and known coefficients; incorporation of this $c_0(x)$ dependence would lead to the {\sc pde} 
\begin{equation}\label{eqn:Westervelt_ssl_nlc_full}
\bigl(\ssl(x) u-\nlc(x)u^2\bigr)_{tt}-\triangle u + \ssl\, D u = \ssl\, r \quad
\mbox{ in }\Omega\times(0,T).
\end{equation}
Neglecting this dependency in the damping term can be justified by smallness of the damping coefficient so that spatial variability of this term has a very minor effect. Neglecting spatial variability of $\ssl$ in the excitation term does not matter due to the fact that the support of $r$ is typically remote from the region of variable (and unknown) sound speed. 

\smallskip

The inverse problem of reconstructing $\nlc(x)$, $\ssl(x)$ from the observations \eqref{eqn:obs} can then be written as 
\begin{equation}\label{eqn:Fnlcsslh}
F(\nlc,\ssl)=h,
\end{equation}
where $F=C\circ S$ and with the parameter-to-state map $S:(\nlc,\ssl)\mapsto u$ where $u$ solves \eqref{eqn:Westervelt_ssl_nlc_intro} and is subject to
the observation operator $C:u\mapsto \mbox{tr}_\Sigma u$.

Well-definedness of the forward operator $F$ and its linearisation in appropriate function spaces is discussed at the beginning of Section~\ref{sec:Newton}.

\subsubsection*{Notation}
Below we will make use of the spaces $\dot{H}^{\beta}(\Omega)$ induced by the norm 
\begin{equation}\label{sobolev_norm}
\|v\|_{H^\beta(\Omega)}=\Bigl(\sum_{j=1}^\infty \lambda_j^\beta\sum_{k\in K^{\lambda_j}}
 |\langle v,\varphi_j^k\rangle|^2\Bigr)^{1/2}
\end{equation}
with the eigensystem $(\lambda_j,\varphi_j)$ of some selfadjoint positive definite operator $\mathcal{A}$ (in this paper, it will be the negative Laplacian with impedance boundary conditions).

Moreover, the Bochner-Sobolev spaces $L^p(0,T;Z)$, $H^q(0,T;Z)$ with $Z$ some Lebesgue or Sobolev spaces and $T$ a finite or infinite time horizon will be used.

We denote the Laplace transform of a function $v\in L^1(0,\infty)$ by $\hat{v}(z)=\int_0^\infty e^{-zt}v(t)\, dt$ for all $z\in\mathbb{C}$ such that this integral exists.

\section{Uniqueness}\label{sec:uniqueness_kappa_c0_lin} 
In this section we will prove linearised uniqueness of $\kappa(x)$ and $c_0(x)$  in $\mathbb{R}^d$ from two observations,  considering the alternative formulation \eqref{eqn:Westervelt_ssl_nlc_intro}, with $T=\infty$ and 
\begin{equation}\label{eqn:def_ssl_nlc}
\ssl=\frac{1}{c_0^2}, \quad \nlc=\frac{\kappa}{c_0^2}.
\end{equation}
To this end we will show injectivity of the linearised forward operator with respect to $\nlc(x)$ and $\ssl(x)$, given two appropriately chosen excitations $\tilde{r}_i$, $i\in\{1,2\}$. On the one hand, this is essential for well-definedness and convergence of the frozen Newton method considered in the reconstruction section below. On the other hand, via \eqref{eqn:def_ssl_nlc}, uniqueness of $\nlc(x)$ and $\ssl(x)$ is equivalent to uniqueness of $\kappa(x)$ and $c_0(x)$.

The linearisation of the forward operator $F:(\nlc,\ssl)\mapsto \mbox{tr}_\Sigma u$ is formally given by 
\hfill\break
$F'(\nlc,\ssl)(\dnlc,\dssl)=\mbox{tr}_\Sigma\du$, where $\du$ solves
\begin{equation}\label{eqn:Westervelt_ssl_nlc_lin}
\begin{aligned}
&\bigl(\ssl(x)\, \du-2\nlc(x)\, u\du\bigr)_{tt}-\triangle \du + \tilde{D} \du = 
-\bigl(\dssl(x)\, u-\dnlc(x)\, u^2\bigr)_{tt} \quad
\mbox{ in }\Omega\times(0,T).
\end{aligned}
\end{equation}

This simplifies considerably if we linearise around vanishing nonlinearity $\nlc=0$ and constant wave speed $\ssl=1/c^2$ for some $c\in\mathbb{R}^+$, which yields $F'(0,\frac{1}{c^2})(\dnlc,\dssl)= \mbox{tr}_\Sigma\du$, where $\du$
solves 
\[
\frac{1}{c^2}\du_{tt}-\triangle \du + \tilde{D} \du = 
-\bigl(\dssl(x)\, u^0-\dnlc(x)\, (u^0)^2\bigr)_{tt}
\mbox{ in }\Omega\times(0,T)
\]
with homogeneous initial and boundary conditions. 
Here $u^0$ solves \eqref{eqn:Westervelt_ssl_nlc_intro} with $\nlc=0$, $\ssl=\frac{1}{c^2}$, which in its turn is a linear constant coefficient \textsc{pde}.

To obtain injectivity of the linearisation, we use two excitations $\tilde{r}_i$, $i\in\{1,2\}$ and the corresponding components of the forward operator $\vec{F}=(F_1,F_2)$ are defined by $F_i=C\circ S_i$ with $S_i:(\nlc,\ssl)\mapsto u_i$ where $u_i$ solves \eqref{eqn:Westervelt_ssl_nlc_intro} with $\tilde{r}=\tilde{r}_i$, $i\in\{1,2\}$ and $C:u\mapsto \mbox{tr}_\Sigma u$.
Our goal is to prove that with an appropriate choice of $\tilde{r}_1$, $\tilde{r}_2$, the only solution to the homogeneous equation $\vec{F}'(0,\frac{1}{c^2})(\dnlc,\dssl)= (0,0)$ is $(\dnlc,\dssl)= (0,0)$. 
To this end, we construct the excitations $\tilde{r}_i$, $i\in\{1,2\}$ such that they lead to space-time separable solutions $u^0_i(x,t)=\phi_i(x)\psi_i(t)$ of \eqref{eqn:Westervelt_ssl_nlc_intro},
\begin{equation}\label{eqn:rtil}
\tilde{r}_i(x,t):= \frac{1}{c^2}\phi(x)\psi_i''(t)-\triangle\phi(x)\psi_i(t) + \tilde{D}[\phi\psi_i](x,t), \quad  i\in\{1,2\}.
\end{equation}

Expanding the solutions $\du_i$, $i\in\{1,2\}$ in terms of eigenfunctions $\varphi_j$ of $-\triangle$,
we can write the Laplace transformed solutions $\du_i$, $i\in\{1,2\}$ as
\begin{equation}\label{eqn:Lapl_du_i}
\Lapl{\underline{du}_i}(x,z)= \sum_{j=1}^\infty \tfrac{1}{\omega_{\lambda_j}}(z) \sum_{k\in K^{\lambda_j}}
\Bigl( \langle \dssl \phi_i, \varphi_j^k\rangle \Lapl{\psi_i''}(z)
+\langle \dnlc \phi_i^2, \varphi_j^k\rangle \Lapl{(\psi_i^2)''}(z) \Bigr) \varphi_j^k(x). 
\end{equation}
Here $(\lambda_j,\varphi_j^k)$ is an eigensystem of $-\triangle$ equipped with the impedance boundary conditions of \eqref{eqn:Westervelt_ssl_nlc_intro}, $\langle\cdot,\cdot\rangle$ is the $L^2$ inner product on $\Omega$,
and 
\[
\omega_\lambda(z)= \begin{cases} 
\frac{1}{c^2}z^2 + \tilde{b}\lambda^\beta z^\alpha + \lambda \mbox{ for \CWKCH}\\
\tilde{b}_2 z^{2+\alpha_2} + \frac{1}{c^2}z^2 + \tilde{b}_1\lambda z^{\alpha_1} + \lambda \mbox{ for \FZ},
\end{cases}
\]
are the reciprocals of the relaxation functions $\frac{1}{\omega_{\lambda}}$.
We will make use of the following two auxiliary results on these relaxation functions.
\begin{lemma}[Lemma 11.4 in \cite{BBB}]\label{lem:what}
For \CWKCH  or \FZ damping, the poles of $\frac{1}{\omega_\lambda}$ differ for different $\lambda$. 
\end{lemma}
\begin{lemma}[Lemma 11.5 in \cite{BBB}]\label{lem:resnon0}
For \CWKCH  or \FZ damping, the residues of the poles of $\frac{1}{\omega_\lambda}$ do no vanish. 
\end{lemma}

With \eqref{eqn:Lapl_du_i}, the premiss $F_i'(0,\frac{1}{c^2})(\dnlc,\dssl)=0$, $i\in\{1,2\}$ reads as 
\begin{equation*}
0= \sum_{j=1}^\infty \tfrac{1}{\omega_{\lambda_j}}(z)
\sum_{k\in K^{\lambda_j}}
\Bigl( \langle \dssl \phi_i, \varphi_j^k\rangle \Lapl{\psi_i''}(z)
+\langle \dnlc \phi_i^2, \varphi_j^k\rangle \Lapl{(\psi_i^2)''}(z)
\Bigr) \varphi_j^k(x_0),\quad x_0\in\Sigma \quad i\in\{1,2\}.
\end{equation*}
Taking the residues at the singularities (which are the poles $p_j$ of the relaxation functions) and applying Lemmas~\ref{lem:what}, \ref{lem:resnon0}, we can single out the contributions pertaining to the individual eigenvalues 
\begin{equation}\label{eqn:Fprime0}
0= 
\sum_{k\in K^{\lambda_j}}
\Bigl( \langle \dssl \phi_i, \varphi_j^k\rangle \Lapl{\psi_i''}(p_j)
+\langle \dnlc \phi_i^2, \varphi_j^k\rangle \Lapl{(\psi_i^2)''}(p_j)
\Bigr) \varphi_j^k(x_0),\quad x_0\in\Sigma \quad i\in\{1,2\},\quad j\in\mathbb{N}.
\end{equation}
In case of one space dimension, the eigenvalues are single and the inner sum consists of one term $\#K^{\lambda_j}=1$. However, in higher space dimensions, we typically have to deal with multidimensional eigenspaces, that is, $\#K^{\lambda_j}>1$, with $(\varphi_j^k)_{k\in K^{\lambda_j}}$ as an orthonormal basis of the eigenspace corresponding to $\lambda_j$.
Looking at each of these eigenspaces individually, it becomes apparent that in order not to lose the essential information separating the individual eigenfunction contributions contained in \eqref{eqn:Fprime0}, we have to make the assumption that these eigenspaces keep their dimension after taking the observation traces.
This can be cast as the linear independence assumption 
\begin{equation}\label{eqn:ass_inj_Sigma_rem}
\left(\sum_{k\in K^\lambda} b_k \varphi_k(x) = 0 \ \mbox{ for all }x\in\Sigma\right)
\ \Longrightarrow \ \left(b_k=0 \mbox{ for all }k\in K^\lambda\right)\,.
\end{equation}
for any eigenvalue $\lambda$ of $-\triangle$ and is basically a geometric condition on $\Sigma$.
Under condition \eqref{eqn:ass_inj_Sigma_rem}, from \eqref{eqn:Fprime0} we immediately obtain
\begin{equation}\label{eqn:Fprime0_contd}
0=  \langle \dssl \phi_i, \varphi_j^k\rangle \Lapl{\psi_i''}(p_j)
+\langle \dnlc \phi_i^2, \varphi_j^k\rangle \Lapl{(\psi_i^2)''}(p_j), 
\quad j\in\mathbb{N}, \quad k\in K^{\lambda_j}, \quad i\in\{1,2\}. 
\end{equation}

Now we set $\phi_1=\phi_2=:\phi$ for some function $\phi\not=0$ almost everywhere in $\Omega$, 
so that for each $k$ and $j$, \eqref{eqn:Fprime0_contd} becomes a two-by-two system of equations for the coefficients $a_j^k:= \langle \dssl \phi, \varphi_j^k\rangle$ and $\,b_j^k:= \langle \dnlc \phi^2, \varphi_j^k\rangle$. 
Choosing $\psi_1$, $\psi_2$ such that for all poles $p_j$, the system matrix is regular, that is, 
\begin{equation}\label{eqn:matrix_j}
0\not=\mbox{det}\left(\begin{array}{cc}
\Lapl{\psi_1''}(p_j)&\Lapl{(\psi_1^2)''}(p_j)\\
\Lapl{\psi_2''}(p_j)&\Lapl{(\psi_2^2)''}(p_j)
\end{array}\right)\quad j\in\mathbb{N},
\end{equation}
we obtain $a_j^k=0$, $b_j^k=0$ for all $j\in\mathbb{N}$, $k\in K^{\lambda_j}$.
Hence, the functions $\dssl\, \phi$ and $\dnlc\,\phi^2$ vanish in $L^2(\Omega)$ and by our choice of $\phi\not=0$ a.e.
this implies that $\dssl=0$, $\dnlc=0$ almost everywhere in $\Omega$.

Thus, we have proven the following.
\begin{theorem}\label{thm:uniqueness_lin_nlc_ssl}
Assume that $T=\infty$, that \eqref{eqn:ass_inj_Sigma_rem} holds for the eigenspaces of $-\triangle$ and that the excitations $\tilde{r}_i$ take the form \eqref{eqn:rtil} with $\phi\in\mathcal{D}(-\triangle)$, $\phi\not=0$ a.e. in $\Omega$ and $\psi_1,\,\psi_2$ satisfying \eqref{eqn:matrix_j}. Then, $F_i'(0,\frac{1}{c^2})(\dnlc,\dssl)=0$, $i\in\{1,2\}$ implies
$\dnlc=0$, $\dssl=0$.
\end{theorem}

The same proof also works with the original $\kappa(x)$ and $c_0(x)$ formulation.

Indeed for $\tilde{F}_i:(\kappa,c_0^2)\mapsto \mbox{tr}_\Sigma u_i$ (note that we take the {\em squared} sound speed as a variable), where $u_i$ solves \eqref{eqn:Westervelt_init_D_intro} with 
\begin{equation}\label{eqn:r_2nd}
r_i(x,t):= \phi(x)\psi_i''(t)-c_0^2\triangle\phi(x)\psi_i(t) + D[\phi\psi_i](x,t), \quad  i\in\{1,2\}
\end{equation}
we get that the linearisation around vanishing nonlinearity coefficient $\kappa(x)=0$ and constant sound speed $c_0^2(x)=c^2$ is $\tilde{F}_i'(0,c^2)(\underline{d\kappa},\underline{d c_0^2})= \mbox{tr}_\Sigma\du_i$, where
\[
\du_{i,tt}-c^2\triangle \du_i + D \du_i = 
\underline{d\kappa}\phi^2(\psi_i^2)''+\underline{d c_0^2}\triangle\phi\psi_i
\ \mbox{ in }\Omega\times(0,T).
\]
Thus, from $\tilde{F}_i'(0,c^2)(\underline{d\kappa},\underline{d c_0^2})=0$ for
$i\in\{1,2\}$, together with \eqref{eqn:ass_inj_Sigma_rem} and
Lemmas~\ref{lem:what}, \ref{lem:resnon0}, we obtain, in place of \eqref{eqn:Fprime0_contd}, that
\[
0=  \langle \underline{d\kappa} \phi, \varphi_j^k\rangle \Lapl{(\psi_i^2)''}(p_j)
+\langle \underline{dc_0^2}\triangle \phi_i^2, \varphi_j^k\rangle \Lapl{\psi_i}(p_j), 
\quad j\in\mathbb{N}, \quad k\in K^{\lambda_j}, \quad i\in\{1,2\}. 
\]
Hence, under the assumption 
\begin{equation}\label{eqn:matrix_j_2nd}
0\not=\mbox{det}\left(\begin{array}{cc}
\Lapl{(\psi_1^2)''}(p_j)&\Lapl{\psi_1}(p_j)\\
\Lapl{(\psi_2^2)''}(p_j)&\Lapl{\psi_2}(p_j)
\end{array}\right)\quad j\in\mathbb{N},
\end{equation}
we obtain the following.

\begin{theorem}\label{thm:uniqueness_lin_kappa_c0}
Assume that $T=\infty$, that \eqref{eqn:ass_inj_Sigma_rem} holds 
and the excitations $r_i$ take the form \eqref{eqn:r_2nd} with $\phi\in\mathcal{D}(-\triangle)$, $\phi\not=0$, $\triangle\phi\not=0$ a.e. in $\Omega$ and $\psi_1,\,\psi_2$ satisfying \eqref{eqn:matrix_j_2nd}. Then, $\tilde{F}_i'(0,c^2)(\underline{d\kappa},\underline{d c_0^2})=0$, $i\in\{1,2\}$ implies
$\underline{d\kappa}=0$, $\underline{d c_0^2}=0$.
\end{theorem}

\section{Reconstruction of the nonlinearity coefficient and sound speed by a regularised Newton scheme}


\subsection{Well-definedeness and convergence a frozen Newton method}\label{sec:Newton}

We first of all restrict ourselves to the classical Kelvin-Voigt damping $\tilde{D}=-\tilde{b}\triangle \partial_t$, that is, \CWKCH with $\alpha=\beta=1$. Later on, in Subsection~\ref{sec:aao}, we will return to both general damping models \CWKCH, \FZ.

By a slight extension of \cite[Theorem 1.1 and Proposition 3]{KL09Westervelt}, the parameter-to-state map 
\begin{equation}\label{eqn:V}
S:\mathcal{D}(F)\to 
V:=H^1(0,T;H^2(\Omega))
\cap W^{1,\infty}(0,T;H^1(\Omega))\cap W^{2,\infty}(0,T;L^2(\Omega))
\end{equation}
is well-defined on 
\[
\mathcal{D}(F):=\{(\nlc,\ssl)\in L^\infty(\Omega)\times L^\infty(0,T;L^\infty(\Omega))\, : \, \frac{1}{\ssl}\in L^\infty(0,T;L^\infty(\Omega)), \ \ssl\in H^1(0,T;L^3(\Omega))\}
\]
for a smooth bounded domain $\Omega\subseteq\mathbb{R}^d$, $d\in\{1,2,3\}$, 
$\tilde{r}\in L^2(0,T;L^2(\Omega))\cup H^1(0,T;H^{-1}(\Omega))$ with $\tilde{r}$ small enough in this norm. 
Note that we will have to deal with a potentially time-dependent $\ssl$ below and thus consider a function space that is able to capture this.
By 
Sobolev's Lemma, this implies that evaluation of $u$ at single points or on a smooth manifold is feasible in a continuous way and thus $F:\mathcal{D}(F)\to Y$ is well-defined for any 
$Y\supseteq L^\infty(0,T;C(\Sigma))$ in case $\Sigma$ is a compact and smooth manifold or $Y\supseteq L^\infty(0,T;\ell^\infty(\Sigma))$ in case $\Sigma$ is a set of discrete points. To make use of a Hilbert space structure, we will simply set $Y=L^2(0,T;L^2(\Sigma))$ or $Y=L^2(0,T;\ell^2(\Sigma))$, respectively.

Likewise it follows that for any $(\nlc,\ssl)\in \mathcal{D}(F)$, $(\dnlc,\dssl)\in L^\infty(\Omega)\times L^\infty(0,T;L^\infty(\Omega))$, the G\^{a}teaux derivative of the forward operator is given by $F'(\nlc,\ssl)(\dnlc,\dssl)=\mbox{tr}_\Sigma\du$, where $\du$ solves \eqref{eqn:Westervelt_ssl_nlc_lin}. 

In particular, for applying a frozen Newton method to \eqref{eqn:Fnlcsslh}, we linearise at $\ssl=1/c^2$ (for some constant $c$), $\nlc=0$, that is, we use $F'(0,1/c^2)(\dnlc,\dssl)=\mbox{tr}_\Sigma\du$, where 
\begin{equation*}
\begin{aligned}
&\frac{1}{c^2}\du_{tt}-\triangle \du + \tilde{D} \du = 
-\bigl(\dssl(x)\, u-\dnlc(x)\, u^2\bigr)_{tt} \quad
\mbox{ in }\Omega\times(0,T)
\end{aligned}
\end{equation*}

Using two well-chosen excitations $\tilde{r}_1$, $\tilde{r}_2$, from Theorem~\ref{thm:uniqueness_lin_nlc_ssl} we have linearised injectivity of the two component forward operator $\vec{F}=(F_1,F_2)$ with  $F_i=C\circ S_i$ and $S_i:(\nlc,\ssl)\mapsto u_i$ defined as the parameter-to-state map for \eqref{eqn:Westervelt_ssl_nlc_intro} with $\tilde{r}=\tilde{r}_i$, $i\in\{1,2\}$. 
Thus we conclude formal well-definedness of a frozen Newton scheme by 
\[
(\nlc_{n+1},\ssl_{n+1})= (\nlc_{n},\ssl_{n}) + (\dnlc,\dssl)\mbox{ where $(\dnlc,\dssl)$ solves }
\vec{F}'(\nlc_0,\ssl_0)(\dnlc,\dssl)= \vec{h}-\vec{F}(\nlc_{n},\ssl_{n})
\]
provided that $\vec{h}-\vec{F}(\nlc_{n},\ssl_{n})$ lies in the range of $\vec{F}'(\nlc_0,\ssl_0)$.
However, the inverse problem inherits the ill-posedness from the original nonlinear one and the given data is typically contaminated with noise, that is, in place of $\vec{h}=(h_1,h_2)$ we only have $\vec{h}^\delta\approx \vec{h}$.
Thus regularisation needs to be applied and the convergence analysis of the resulting iterative reconstruction scheme requires structural conditions on the forward operator. 
One of the conditions allowing for convergence guarantees is range invariance of the linearised forward operator (as plausible from the above requirement of the residual lying in the range of $\vec{F}'(\nlc_0,\ssl_0)$)
and can be established for our problem in the slightly relaxed form 
\begin{equation}\label{rangeinvar_diff_vec}
\vec{F}(\nlc,\vec{\ssl})-\vec{F}(\nlc_0,\vec{\ssl}_0)=\vec{F}'(\nlc_0,\vec{\ssl}_0)(\dnlc(\nlc,\vec{\ssl}),\vec{\dssl}(\nlc,\vec{\ssl})).
\end{equation}
To illustrate this first of all in the single excitation case 
\begin{equation}\label{rangeinvar_diff}
F(\nlc,\ssl)-F(\nlc_0,\ssl_0)=F'(\nlc_0,\ssl_0)(\dnlc(\nlc,\ssl),\dssl(\nlc,\ssl)),
\end{equation}
note that it is straightforward to see that by setting 
\begin{equation}\label{rangeinvar_defr}
\begin{aligned}
\dnlc(\nlc,\ssl)&=\dnlc(\nlc):=\nlc-\nlc_0, \\
\dssl(\nlc,\ssl)&=\dssl(\nlc,\ssl;u,u_0):=\frac{1}{u_0}\Bigl((\ssl-\ssl_0)u-(\nlc-\nlc_0)(u^2-u_0^2)-\nlc_0(u-u_0)^2\Bigr)
\end{aligned}
\end{equation}
we can satisfy the identity \eqref{rangeinvar_diff}.
However, through time dependence of $u_0$ and $u$, the expression for $\dssl(\nlc,\ssl)$ in the second identity of \eqref{rangeinvar_defr} will be time dependent as well.
Thus we consider $\ssl$ as a space {\em and} time dependent function. 
Moreover, we have to take into account two excitations resulting in two different states $u^i$, $i\in\{1,2\}$ and thus also need two copies of $\ssl$ to be able to capture this in \eqref{rangeinvar_defr}, thus considering $\vec{\ssl}(x,t)=(\ssl^1(x,t),\ssl^2(x,t))$.
Introducing so much additional dimensionality in parameter space clearly counteracts uniqueness and after all our aim is to reconstruct only one $\ssl(x)$ depending only on the space variables (along with the nonlinearity coefficient $\nlc(x)$).
This is achieved by penalisation with an operator 
\begin{equation}\label{eqn:P}
P:(\ssl^1,\ssl^2)\mapsto 
(\ssl^1 - \text{Proj}^{L^2(0,T;\mu)}_{\textup{const}}\ssl^1,
\ssl^2 - \text{Proj}^{L^2(0,T;\mu)}_{\textup{const}}\ssl^2) 
\end{equation}
where $\text{Proj}^{L^2_\mu(0,T)}_{\textup{const}}$ is the $L^2_\mu$ projection on the space of constant functions with a finite measure $\mu$ on $(0,T)$, including the case of an infinite time horizon $T=\infty$.
Note that in the latter case we do not use $\mu$ as the ordinary Lebesgue measure $\lambda$, since this would exclude the constant-in-time solutions that we are actually looking for but, e.g., define $\mu$ by $d\mu(t)=t^{-2}d\lambda(t)$.   

Setting $(\nlc_0,\vec{\ssl}_0)=(0,1/c^2,1/c^2)$ and abbreviating
\[
\begin{aligned}
&K:=\vec{F}'(\nlc_0,\vec{\ssl}_0), \quad \vec{h}_0=\vec{h}-\vec{F}(\nlc_0,\vec{\ssl}_0), \\
&r:(\nlc,\ssl^1,\ssl^2)\mapsto (\dnlc(\nlc),\dssl(\nlc,\ssl^1;u^1,u^1_0),\dssl(\nlc,\ssl^2;u^2,u^2_0))
\mbox{ as in \eqref{rangeinvar_defr}},
\end{aligned}
\]
we can thus write the original inverse problem \eqref{eqn:Fnlcsslh} equivalently as a combination of an ill-posed linear and a well-posed nonlinear problem   
\begin{equation}\label{FP}
\begin{aligned}
&K\hat{r}=\vec{h}_0\\
&r(\nlc,\vec{\ssl})=\hat{r}\\
&P\vec{\ssl}=0
\end{aligned}
\end{equation}
for the unknowns $(\nlc,\vec{\ssl},\hat{r})$.
In view of \eqref{rangeinvar_defr} we expect $r$ to be close to the identity in the sense of 
\begin{equation}\label{rid}
\exists\, C_r\in(0,1)\, \forall (\nlc,\vec{\ssl})\in U(\subseteq X)\, : \;
\|r(\nlc,\vec{\ssl})-(\nlc-\nlc_0,\vec{\ssl}-\vec{\ssl}_0)\|_X
\rerevision{\leq C_r\|(\nlc-\nlc_0,\vec{\ssl}-\vec{\ssl}_0)\|_X^2}
\end{equation}
\Margin{Ref 2}
for some sufficiently small neighborhood $U\subseteq X$ of the exact solution $(\nlc^\dagger,\vec{\ssl}^{\ \dagger})$, an estimate that we will establish in an appropriate function space setting $X$ in the proof of Theorem~\ref{thm:convfrozenNewton} below.

Thus, a natural way of making use of the structure \eqref{FP} in a regularised frozen Newton type method is to define iterates for $\vec{x}=(\nlc,\ssl^1,\ssl^2)$ with $r(\rerevision{\vec{x}})\approx \rerevision{\vec{x}}-\rerevision{\vec{x}}_0$ as minimizers of
\Margin{Ref 1 1.}
\begin{equation}\label{frozenNewton}
\rerevision{\vec{x}}_{n+1}^\delta \in \mbox{argmin}_{\rerevision{\vec{x}}\in \mathcal{D}(\vec{F})}
\|K(\rerevision{\vec{x}}-\rerevision{\vec{x}}_n^\delta)+\vec{F}(\rerevision{\vec{x}}_n^\delta)-h^\delta\|_Y^p+\regpar_n\mathcal{R}(\rerevision{\vec{x}})+\|P\rerevision{\vec{x}}\|_Z^2.
\end{equation}
with a proper $X$ lower semicontinuous functional $\mathcal{R}$, a sequence of positive regularisation parameters tending to zero 
\Margin{Ref 1 3.}
\rerevision{$(\regpar_n)_{n\in\mathbb{N}}\subseteq\mathbb{R}^+$, $\regpar_n\stackrel{n\to\infty}{\longrightarrow}0$ \footnote{not to be mistaken with the boundary impedance coefficient in \eqref{eqn:Westervelt_init_D_intro}}} and $P$ as in \eqref{eqn:P}, $Z:=L^2_\mu(0,T;L^2(\Omega))^2$.
In view of the well-posedness results quoted above, we choose $X\subseteq L^\infty(\Omega)\times L^\infty(0,T;L^\infty(\Omega))^2$ and $Y\supseteq L^\infty(0,T;L^\infty_\nu(\Sigma))$ with $\nu$ being just the Lebesgue measure in case of a smooth manifold $\Sigma$ and the counting measure in case $\Sigma$ consists of discrete points.

This includes the Hilbert space setting 
\begin{equation}\label{eqn:X}
\begin{aligned}
&X = X_\nlc\times X_\ssl^2 \textup{ with }
X_\nlc=H^\sigma(\Omega), \ X_\ssl= H^\tau(0,T;H^1(\Omega))\cap L^2(0,T;H^2(\Omega)),  \\ 
&Y = L^2(0,T;L^2_\nu(\Sigma))^2
\end{aligned}
\end{equation}
in space dimensions $d\in\{1,2,3\}$ for any $\sigma>d/2$, $\tau>1/2$,  
since due to Sobolev's and Agmon's interpolation inequality 
\footnote{$\|v\|_{L^\infty(\Omega)}^2\leq C\|v\|_{H^1(\Omega)}\|v\|_{H^2(\Omega)}$}, 
$X_\ssl$ continuously embeds into $L^\infty(\Omega)\times L^\infty(0,T;L^\infty(\Omega))$.
With this and $\mathcal{R}(\vec{x})=\|\vec{x}-\vec{x}_0\|_X^2$ we can write \eqref{frozenNewton} in terms of its necessary (and due to convexity, sufficient) first order optimality conditions
\begin{equation}\label{frozenNewtonHilbert}
\vec{x}_{n+1}^\delta=\vec{x}_{n}^\delta+(K^\star K+P^\star P+\regpar_n I)^{-1}
\Bigl(K^\star (\vec{h}^\delta-\vec{F}(\vec{x}_n^\delta))-P^\star P\vec{x}_n^\delta+\regpar_n(\vec{x}_0-\vec{x}_n^\delta)\Bigr)
\end{equation}
where $K^\star $ denotes the Hilbert space adjoint of $K:X\to Y$. 

In case of noisy data with noise level $\delta$ satisfying the bound 
\begin{equation}\label{eqn:delta}
\|\vec{h}^\delta-\vec{F}(\nlc^\dagger,\ssl^\dagger)\|_{L^2(0,T;L^2(\Omega))}\leq\delta
\end{equation}
we have to stop the iteration according to 
\begin{equation}\label{nstar}
n_*(\delta)\to0, \quad \delta\sum_{j=0}^{n_*(\delta)-1}c_r^j\regpar_{n_*(\delta)-j-1}^{-1/2} \to 0 \qquad \textup{ as }\delta\to0
\end{equation}
with $c_r$ as in \eqref{rid}.
With the simple geometric sequence $\regpar_n=\regpar_0\theta^n$ for some 
$\theta\in(0,1)$, 
this just corresponds to the usual a priori choice  $\regpar_{n_*(\delta)}\to0$ and $\delta^2/\regpar_{n_*(\delta)-1}\rerevision{\to0}$ \rerevision{ as $\delta\to0$, cf., e.g., \cite[Chapter 4]{KalNeuSch08}}.
\Margin{Ref 1 2.}

In view of \eqref{rangeinvar_defr}, we also 
assume that $u_0^i=S^i(\nlc_0,\ssl_0^i)$ is bounded away from zero
\begin{equation}\label{u0}
|u_0^i|\geq \underline{c}>0\ \mbox{ a.e. on }(0,T)\times\Omega, \ i\in\{1,2\},
\end{equation}
We will comment on this condition in Section~\ref{sec:aao} below.

From \cite[Theorem 2.2]{rangeinvar} and Theorem~\ref{thm:uniqueness_lin_nlc_ssl} we thus conclude the following convergence result.

\begin{theorem}\label{thm:convfrozenNewton}
Let the conditions of Theorem~\ref{thm:uniqueness_lin_nlc_ssl} on the observation set $\Sigma$ and on the excitations $\tilde{r}^1,\, \tilde{r}^2$ be satisfied.
Let $\vec{x}_0=(\nlc_0,\ssl_0^1,\ssl_0^2)
\in U:=\mathcal{B}_\rho(\vec{x}^\dagger)$ for some $\rho>0$ sufficiently small, assume that \eqref{u0} holds for $u_0^i=S^i(\nlc_0,\ssl_0^i)$, $i\in\{1,2\}$ 
and let the stopping index $n_*=n_*(\delta)$ be chosen according to \eqref{nstar}.

Then the iterates $(\vec{x}_n^\delta)_{n\in\{1,\ldots,n_*(\delta)\}}$ are well-defined by \eqref{frozenNewtonHilbert}, remain in $\mathcal{B}_\rho(\vec{x}^\dagger)$ and converge in $X$ 
(defined as in \eqref{eqn:X} with $\tau\in(1,5/4)$), 
$\|\vec{x}_{n_*(\delta)}^\delta-\vec{x}^\dagger\|_X\to0$ as $\delta\to0$. In the noise free case $\delta=0$, $n_*(\delta)=\infty$ we have $\|\vec{x}_n-\vec{x}^\dagger\|_X\to0$ as $n\to\infty$.
\end{theorem}

{\em Proof. \ }
With $\tau\in(1,5/4)$, the solution space $V$ defined in \eqref{eqn:V} is
embedded in $H^\tau(0,T;L^\infty(\Omega)\cap W^{1,3}(\Omega))$ and on the other hand the parameter space $X_\ssl$ defined in \eqref{eqn:X} is embedded in $L^\infty(0,T;L^\infty(\Omega))\cap H^\tau(0,T;L^6(\Omega))$; moreover $\nlc-\nlc_0$ and $\nlc_0$ are time-independent.
Using the fact that $H^\tau(0,T)$ and $H^2(\Omega)$ are Banach algebras, we obtain the following estimates.
We start from the identity
\[
\dssl(\nlc,\ssl)-(\ssl-\ssl_0)
= \tfrac{1}{u_0}(u-u_0)\Bigl((\ssl-\ssl_0)-(\nlc-\nlc_0)(u+u_0)-\nlc_0(u-u_0)\Bigr),
\]
and, using the above mentioned continuous embeddings, estimate the individual terms as follows
\[
\begin{aligned}
&\|\tfrac{1}{u_0}(u-u_0)(\ssl-\ssl_0)\|_{L^2(H^2)}\leq
\|\tfrac{1}{u_0}\|_{L^\infty(H^2)} \|u-u_0\|_{L^\infty(H^2)} \|\ssl-\ssl_0\|_{L^2(H^2)}\\
&\|\tfrac{1}{u_0}(u-u_0)(\ssl-\ssl_0)\|_{H^\tau(H^1)}\leq
\|\tfrac{1}{u_0}\|_{H^\tau(L^\infty)} \|u-u_0\|_{H^\tau(L^\infty)} \|\ssl-\ssl_0\|_{H^\tau(H^1)}\\
&\quad+\Bigl(\|\tfrac{1}{u_0}\|_{H^\tau(L^\infty)} \|(u-u_0)\|_{H^\tau(W^{1,3})} 
+\|\tfrac{1}{u_0}\|_{H^\tau(W^{1,3})} \|u-u_0\|_{H^\tau(L^\infty)}\Bigr) \|\ssl-\ssl_0\|_{H^\tau(L^6)}
\end{aligned}
\]
\[
\begin{aligned}
&\|\tfrac{1}{u_0}(u-u_0)(u+u_0)(\nlc-\nlc_0)\|_{L^2(H^2)}\leq
\|\tfrac{1}{u_0}\|_{L^2(H^2)} \|u-u_0\|_{L^\infty(H^2)} \|u+u_0\|_{L^\infty(H^2)} \|\nlc-\nlc_0\|_{H^2}\\
&\|\tfrac{1}{u_0}(u-u_0)(u+u_0)(\nlc-\nlc_0)\|_{H^\tau(H^1)}\leq
\|\tfrac{1}{u_0}\|_{H^\tau(L^\infty)} \|u-u_0\|_{H^\tau(L^\infty)} \|u+u_0\|_{H^\tau(L^\infty)} \|\nlc-\nlc_0\|_{H^1}\\
&\hspace*{4cm}+\Bigl(\|\tfrac{1}{u_0}\|_{H^\tau(L^\infty)} \|u-u_0\|_{H^\tau(L^\infty)} \|(u+u_0)\|_{H^\tau(W^{1,3})} \\
&\hspace*{4.5cm}+\|\tfrac{1}{u_0}\|_{H^\tau(L^\infty)} \|u-u_0\|_{H^\tau(W^{1,3})} \|(u+u_0)\|_{H^\tau(L^\infty)}\\
&\hspace*{4.5cm}+\|\tfrac{1}{u_0}\|_{H^\tau(W^{1,3})} \|u-u_0\|_{H^\tau(L^\infty)} \|(u+u_0)\|_{H^\tau(L^\infty)}\Bigr) \|\nlc-\nlc_0\|_{H^\tau(L^6)}
\end{aligned}
\]
and analogously for the term $\tfrac{1}{u_0}(u-u_0)^2\nlc_0$.

This together with assumption \eqref{u0} and Lipschitz continuity of $S$ allow us to establish
\[
\|(\dnlc(\nlc,\vec{\ssl}),\vec{\dssl}(\nlc,\vec{\ssl}))-(\nlc-\nlc_0,\vec{\ssl}-\vec{\ssl}_0)\|_X\leq C_r\|(\nlc-\nlc_0,\vec{\ssl}-\vec{\ssl}_0)\|_X^2
\]
\Margin{Ref 2}
\rerevision{that is, \eqref{rid}, which together with $r(\nlc_0,\vec{\ssl}_0)=0$ implies that $r$ is Fr\'{e}chet differentiable at $(\nlc_0,\vec{\ssl}_0)$ with derivative $r(\nlc_0,\vec{\ssl}_0)=\textup{id}_X$. 
Applying analogous estimates 
\arxiv{ to the Taylor remainder decomposition $r(\nlc+\delnlc,\vec{\ssl}+\vec{\delssl}) - r(\nlc,\vec{\ssl})-r'(\nlc,\vec{\ssl})(\delnlc,\vec{\delssl})=(0,d_{\textup{Tay}1},d_{\textup{Tay}2})$ with (skipping the indices again)
\[
\begin{aligned}
d_{\textup{Tay}}=\tfrac{1}{u_0}\Bigl(&
(\ssl-\ssl_0)(\tilde{u}-u-\delu)+\delssl(\tilde{u}-u)\\
&-(\nlc-\nlc_0)(\tilde{u}+u)(\tilde{u}-u-\delu)+\delu(\tilde{u}-u)-\delnlc(\tilde{u}+u)(\tilde{u}-u)\\
&-\nlc_0\bigl((\tilde{u}+u-2u_0)(\tilde{u}-u-\delu)+\delu(\tilde{u}-u)\bigr)\Bigr),
\end{aligned}
\]
where $\tilde{u}=S(\nlc+\delnlc,\ssl+\delssl)$, $u=S(\nlc,\ssl)$, 
$\delu=S'(\nlc,\ssl)(\delnlc,\delssl)$,}
one can establish continuous differentiability of $r$ and thus, with 
\[
\begin{aligned}
&\|r(\nlc,\vec{\ssl})-r(\nlc^\dagger,\vec{\ssl}^\dagger)-(\nlc-\nlc^\dagger,\vec{\ssl}-\vec{\ssl}^\dagger)\|_X
=\|r(\nlc,\vec{\ssl})-r(\nlc^\dagger,\vec{\ssl}^\dagger)-r'(\nlc_0,\vec{\ssl}_0)(\nlc-\nlc^\dagger,\vec{\ssl}-\vec{\ssl}^\dagger)\|_X\\
&\leq\|r(\nlc,\vec{\ssl})-r(\nlc^\dagger,\vec{\ssl}^\dagger)-r'(\nlc^\dagger,\vec{\ssl}^\dagger)(\nlc-\nlc^\dagger,\vec{\ssl}-\vec{\ssl}^\dagger)\|_X
+\|[r'(\nlc^\dagger,\vec{\ssl}^\dagger)-r'(\nlc_0,\vec{\ssl}_0)](\nlc-\nlc^\dagger,\vec{\ssl}-\vec{\ssl}^\dagger)\|_X
\end{aligned}
\]
a bound of the form
\begin{equation}\label{rid_cr}
\begin{aligned}
\exists\, c_r\in(0,1)\, \forall (\nlc,\vec{\ssl})\in U(\subseteq X)\, : \;&
\|r(\nlc,\vec{\ssl})-r(\nlc^\dagger,\vec{\ssl}^\dagger)-(\nlc-\nlc^\dagger,\vec{\ssl}-\vec{\ssl}^\dagger)\|_X\\
&\leq c_r\|(\nlc-\nlc^\dagger,\vec{\ssl}-\vec{\ssl}^\dagger)\|_X
\end{aligned}
\end{equation}
on a sufficiently small neighborhood $U$ of $(\nlc^\dagger,\vec{\ssl}^\dagger)$ containing $(\nlc_0,\vec{\ssl}_0)$, 
which is the crucial convergence condition in \cite[Theorem 2.2]{rangeinvar}.
} 
\begin{flushright}
$\diamondsuit$
\end{flushright}

\subsubsection{An all-at-once formulation}\label{sec:aao}
As can be seen from the proof of Theorem~\ref{thm:convfrozenNewton}, we need to avoid division by zero by assuming \eqref{u0}; however, as a solution to a wave equation, $u_0$ will typically change sign. 
This problem can be circumvented by considering the all-at-once version, which
allows us to choose $u_0$ not necessarily as a \textsc{pde} solution and also provides 
for much more freedom in the choice of the function spaces.

To this end, we consider the model and the observation equation as a system of operator equations for the sought-after coefficients and the states. That is, we set $\vec{x}=(\nlc,\ssl^1,\ssl^2,u^1,u^2)$ and replace the definition of $\vec{F}=(F_1,F_2)$ by 
\[
F_i(\nlc,\ssl^1,\ssl^2,u^1,u^2) = 
\left(\begin{array}{l}F^{\textup{mod}}(\nlc,\ssl^i,u^i)\\ \textup{tr}_\Sigma u^i\end{array}\right) 
\]
where 
\begin{equation}\label{eqn:veryweakmodel}
\begin{aligned}
&\langle F^{\textup{mod}}(\nlc,\ssl,u), w\rangle_{W^*,W}
= \int_0^T \int_\Omega\Bigl(\bigl(\ssl(x) u-\nlc(x)u^2\bigr)\, w_{tt}+(-\triangle u + \tilde{D} u)\,w\Bigr)\, dx\, dt
\\
&\quad w\in W:= \{v\in H^2(0,T;L^2(\Omega))\ : \ v(T)=0, \ v_t(T)=0\}, 
\end{aligned}
\end{equation}
where $\langle\cdot,\cdot\rangle_{W^*,W}$ denotes the dual pairing in $W$.\relax
\footnote{for $T=\infty$ the end conditions are to be understood in a limiting sense $\lim_{t\to\infty} v(t)=0$, $\lim_{t\to\infty} v_t(t)=0$}
Note that we aim here to allow for low regularity of the coefficients to decrease the degree of ill-posedness of the inverse problem as much as possible. This is enabled by the weak formulation (with respect to time derivatives) of the \textsc{pde} model in \eqref{eqn:veryweakmodel} and allows to use the function spaces 
\begin{equation}\label{eqn:X_aao}
\begin{aligned}
&X = X_\nlc\times X_\ssl^2\times X_u^2 \textup{ with }
X_\nlc=L^2(\Omega), \ X_\ssl= L^2(0,T;L^2(\Omega)), \\ 
&X_u= \begin{cases}L^2(0,T;H^2(\Omega))\cap H^{\tilde{\alpha}}(0,T;\dot{H}^{\tilde{\beta}}(\Omega)) \mbox{ in case of \CWKCH }\\
H^{\tilde{\alpha_1}}(0,T;H^2(\Omega))\cap H^{\tilde{\alpha_1}+2}(0,T;L^2(\Omega))  \mbox{ in case of \FZ }\end{cases}\\
&Y = (W^*)^2\times L^2(0,T;L^2_\nu(\Sigma))^2.
\end{aligned} 
\end{equation}
Here the exponents are chosen such that $\tilde{D}$ maps $X_u$ into $L^2(0,T;L^2(\Omega))$ and 
$X_u$ continuously embeds into $L^\infty(0,T;L^\infty(\Omega))$:
\begin{equation}\label{eqn:alphabetatilde}
\tilde{\alpha}\geq\alpha, \ 
\tilde{\beta}\geq\beta, \
\tilde{\alpha}_1\geq\alpha_1, \
\tilde{\alpha}_2\geq\alpha_2, \
\tilde{\alpha}>\frac12, \ \tilde{\beta}>\frac{d}{2}, \
\max\{\tilde{\alpha}_1,\tilde{\alpha}_2\}>0\,.
\end{equation}

The range invariance condition
\begin{equation*}
\vec{F}(\nlc,\vec{\ssl},\vec{u})-\vec{F}(\nlc_0,\vec{\ssl}_0,\vec{u}_0)=\vec{F}'(\nlc_0,\vec{\ssl}_0,\vec{u}_0)(\dnlc(\nlc,\vec{\ssl}),\vec{\dssl}(\nlc,\vec{\ssl})\vec{\du}(\nlc,\vec{\ssl}))
\end{equation*}
can be verified with $\dnlc$, $\dssl$ defined as in \eqref{rangeinvar_defr} and $\du=u-u_0$.
Injectivity of $\vec{F}'(0,1/c^2,1/c^2,u_0^1,u_0^2)$ can be shown analogously to Theorem~\ref{thm:uniqueness_lin_nlc_ssl}, 
for $u^0_i(x,t)=\phi(x)\psi_i(t)$ with 
$\phi\in\mathcal{D}(-\triangle)$, $\phi\not=0$ a.e. in $\Omega$, $\psi_1,\,\psi_2$ satisfying \eqref{eqn:matrix_j},
but actually without $u^0_i$ needing to solve \eqref{eqn:Westervelt_ssl_nlc_intro}.
Finally, the estimate of $r(\vec{x})-(\vec{x}-\vec{x}_0)$ simplifies to
\[
\begin{aligned}
&\|r(\vec{x})-(\vec{x}-\vec{x}_0)\|_X^2 
= \sum_{i=1}^2\|\tfrac{1}{u_0}(u-u_0)\Bigl((\ssl-\ssl_0)-(\nlc-\nlc_0)(u+u_0)-\nlc_0(u-u_0)\Bigr)\|_{L^2(L^2)}^2\\
&\leq \sum_{i=1}^2\|\tfrac{1}{u_0}\|_{L^\infty(L^\infty)}^2 \|(u-u_0)\|_{L^\infty(L^\infty)}^2 
\Bigl(\|\ssl-\ssl_0\|_{L^2} \\
&\hspace*{3cm}+ \|\nlc-\nlc_0\|_{L^2(L^2)} \|u+u_0\|_{L^\infty(L^\infty)}
+\|\nlc_0\|_{L^2(L^2)} \|u-u_0\|_{L^\infty(L^\infty)} \Bigr)^2.
\end{aligned}
\]

Thus, applicability and convergence of the frozen Newton method transfers to this all-at-once setting as follows.

\begin{theorem}\label{thm:convfrozenNewton_aao}
Let the conditions of Theorem~\ref{thm:uniqueness_lin_nlc_ssl} on the observation set $\Sigma$ and on the functions $\phi$, $\psi_1$, $\psi_2$ in $u^0_i(x,t)=\phi(x)\psi_i(t)$ be satisfied.
Let $\vec{x}_0=(\nlc_0,\ssl_0^1,\ssl_0^2,u_0^1,u_0^2)\in U:=\mathcal{B}_\rho(\vec{x}^\dagger)$ for some $\rho>0$ sufficiently small, assume that \eqref{u0} holds
and let the stopping index $n_*=n_*(\delta)$ be chosen according to \eqref{nstar}.

Then the iterates $(\vec{x}_n^\delta)_{n\in\{1,\ldots,n_*(\delta)\}}$ are well-defined by \eqref{frozenNewtonHilbert}, remain in $\mathcal{B}_\rho(\vec{x}^\dagger)$ and converge in $X$ 
(defined as in \eqref{eqn:X_aao} with \eqref{eqn:alphabetatilde}), 
$\|\vec{x}_{n_*(\delta)}^\delta-\vec{x}^\dagger\|_X\to0$ as $\delta\to0$.
In the noise-free case $\delta=0$, $n_*(\delta)=\infty$
we have $\|\vec{x}_n-\vec{x}^\dagger\|_X\to0$ as $n\to\infty$.
\end{theorem}

The price to pay for this more relaxed setting is convergence of the
coefficients in a weaker norm as compared to Theorem~\ref{thm:convfrozenNewton}.

\subsection{Reconstructions}\label{sec:reconstructions}

In this section we show reconstructions of $\nlc(x)$, $\ssl(x)$ in \eqref{eqn:Westervelt_ssl_nlc_intro} with Caputo-Wismer-Kelvin damping, that is, 
\begin{equation}\label{eqn:Westervelt_ssl_nlc_CH}
\begin{aligned}
&\bigl(\ssl(x) u-\nlc(x)u^2\bigr)_{tt}-\triangle u - b\triangle \partial_t^\alpha u = \tilde{r} \quad
\mbox{ in }\Omega\times(0,T)\\
\partial_\nu u+\gamma u&=0 \mbox{ on }\partial\Omega\times(0,T),\quad
u(0)=0, \quad u_t(0)=0 \quad \mbox{ in }\;\Omega.
\end{aligned}
\end{equation}
in one space dimension $\Omega=(0,1)$ with Dirichlet-Neumann conditions $\gamma(0)=\infty$, $\gamma(1)=0$ from measurements at two points $\Sigma=(0.1,1)$. (Note that since we impose homogeneous Dirichlet boundary conditions at the left endpoint, measuring $u$ there would not provide any additional information; indeed, also in practice the transducer array will be immersed into the overall computational domain $\Omega$.)

For the numerical solution of \eqref{eqn:Westervelt_ssl_nlc_CH},
we follow the method of \cite{nonlinearity_imaging_fracWest} and rewrite the
equation by integrating once with respect to time
\begin{equation}\label{eqn:Westervelt_ssl_nlc_int}
\begin{aligned}
&\bigl(\ssl(x) -2\nlc(x)u\bigr)u_{t}-\tilde{b}\triangle I_t^{1-\alpha}u - \triangle I_t^{1}u = I_t^{1}\tilde{r} \quad
\mbox{ in }\Omega\times(0,T)\\
\partial_\nu u+\gamma u&=0 \mbox{ on }\partial\Omega\times(0,T),\quad
u(0)=0 \quad \mbox{ in }\;\Omega,
\end{aligned}
\end{equation}
to which we apply a modified Crank-Nicolson solver taking into account the fractional integral term. 
Likewise for its linearisation \eqref{eqn:Westervelt_ssl_nlc_lin}.

\medskip

To test the \rerevision{(reduced) frozen Newton method \eqref{frozenNewton}} 
\Margin{Ref 1 4.}
from Section~\ref{sec:Newton}, we consider three scenarios 
\rerevision{(a), (b) and (c)} 
\Margin{Ref 1 5.}
as described below.
While the theory from Section~\ref{sec:Newton} requires two excitations
and an extension  of $\ssl$ to a time dependent function,
this was not needed in practical computations. The reconstructions shown here are based on a single excitation, but carrying out measurements at two points $\Sigma=\{0.1,1\}$. Also $\ssl$ is treated as a function of $x$ only.  
Both coefficients were discretised using a chapeau basis set and the starting values were $\nlc_0=0$ and $\ssl=1$.

\begin{figure}[h]
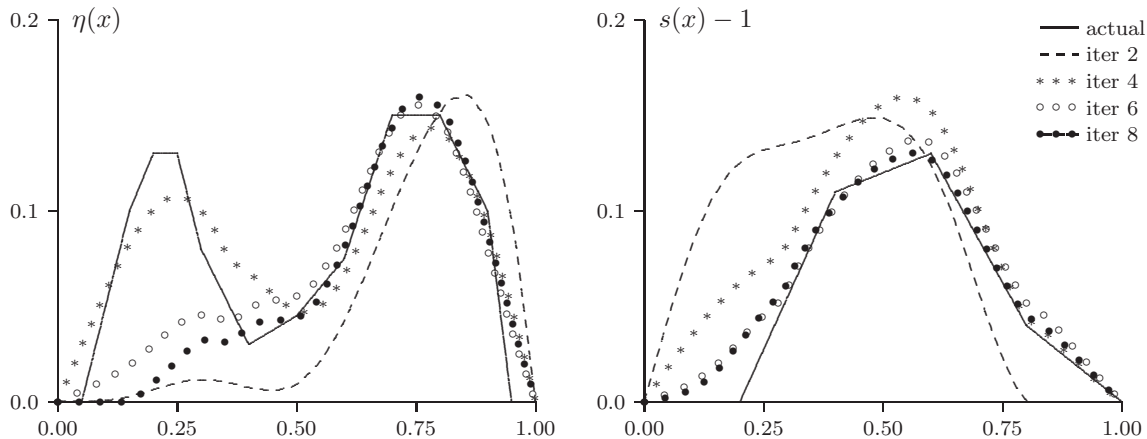

\hbox to \hsize{\hss\copy\figureseven\hss\copy\figureeight\hss}
\centering
\caption{\small {\bf Reconstruction of both $\nlc(x)$ and $\ssl(x)$
under $0.1\%$ data noise; test case (c).} }
\label{fig:kappa+slow_noise01_C}
\end{figure}

Figures \ref{fig:kappa+slow_noise1} and \ref{fig:kappa+slow_noise01}
show a simultaneous reconstruction of both $\nlc(x)$ and $\ssl(x)$ under
$1\%$ and $0.1\%$ noise in the time trace data.
Here the value of the solution $u(x,t)$  was negative and there was therefore
no cancellation effect on the  $\ssl(x)$ and the $\nlc(x)$ term (test case (a)).
\begin{figure}[h]
\vspace*{0.8cm}
\hbox to \hsize{\hss\copy\figureone\hss\copy\figuretwo\hss}
\centering
\caption{\small {\bf Reconstruction of both $\nlc(x)$ and $\ssl(x)$
under $1\%$ data noise; test case (a).}}
\label{fig:kappa+slow_noise1}
\end{figure}
\begin{figure}[h]
\vspace*{0.8cm}
\hbox to \hsize{\hss\copy\figurethree\hss\copy\figurefour\hss}
\centering
\caption{\small {\bf Reconstruction of both $\nlc(x)$ and $\ssl(x)$ under $0.1\%$ data noise; test case (a).}
\label{fig:kappa+slow_noise01}}
\end{figure}
In Figure~\ref{fig:kappa+slow_noise01_B} we show the difference
when the sign of $u(x,t)$ is reversed so that there is the potential
for a cancellation effect between $\nlc(x)$ and $\ssl(x)$ (test case (b)).
In fact this occurred resulting in a poorer reconstruction in both
functions.
 Data noise here was again $0.1\%$.
\begin{figure}[h]
\hbox to \hsize{\hss\copy\figurefive\hss\copy\figuresix\hss}
\centering
\caption{\small {\bf Reconstruction of both $\nlc(x)$ and $\ssl(x)$
under $0.1\%$ data noise; test case (b).} }
\label{fig:kappa+slow_noise01_B}
\end{figure}
The final picture \ref{fig:kappa+slow_noise01_C}
shows a more complex function $\nlc(x)$ with two features (test case (c)),
one at each end of the interval.
For this run the function $u(x,t)$ was zero at the endpoint $x=0$ and
so small in comparison at the left end as opposed to the right.
Since $\nlc$ occurs in combination with $u$ in the equation this
means a relative loss of information at the left-hand endpoint.
This is clearly visible from the left hand graphic.
In, addition this error in $\nlc(x)$ now affects the
combined term $(\ssl - \nlc u)$ and results in a similarly poor
reconstruction of $\ssl(x)$ near $x=0$.
Note that a seemingly overall better match of $\nlc$ at the fourth iterate is dismissed in subsequent iterations that are much worse in approximating the left hand feature.
This is due to the fact that the mismatch is weighed by the values of $u$ which are small near the left endpoint, but penalize deviations occurring in the right half of the interval (as is the case for iteration 4) much more strongly.
These reconstructions were made under $0.1\%$ data noise 
\rerevision{ with $\regpar_n=0.8^n\cdot 10^{-3}$}
and the discrepancy principle was used as a stopping criterion 
\rerevision{ which basically gave the same results as the a priori choice \eqref{nstar}}. 
\Margin{Ref 1 3.}

This effect of smallness of $u$ at the left hand endpoint is also
apparent in the other figures, particularly in the case where there is
a significant feature near this endpoint.
In all figures we imposed the sign constraint imposed by the physical
problem that $\nlc(x)\geq 0$.

For an $\nlc$ function with support away from $x=0$ the reconstructions
shown in figure \ref{fig:kappa+slow_noise01} under $0.1\%$ noise and in
figure \ref{fig:kappa+slow_noise1} with $1\%$ added noise
indicate a reasonable reconstruction of both $\nlc$ and $\ssl(x)$.
Note that a poor initial guess (both these functions taken to be constant zero and one respectively)
leads to a severe overshoot in the first computed approximation to $\ssl(x)$
although this quickly settles down.
In this case both actual $\ssl(x)$ and $\nlc(x)$ functions have
support away from the left-hand endpoint $x=0$ and there is also overshoot in
the first iteration of $\nlc$.

The difference between $(0.1\%)$ and $1\%$ of added noise
to the data simulated by the direct solver is quite apparent and indicates
the severe ill-conditioning of the inverse problem.

\begin{wrapfigure}{r}{0.45\textwidth}
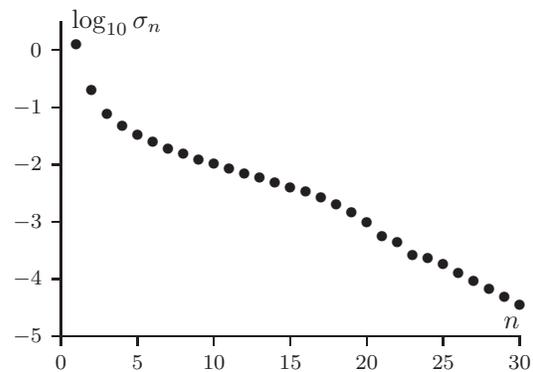

\vspace{-6pt}
\hbox to \hsize{\hss\box\figurenine\hss}
\caption{\small {\bf Singular values of the Jacobian.}}
\label{fig:sv_jacobian}
\end{wrapfigure}
Finally, we show a plot of the singular values of the
Jacobian matrix frozen at $\ssl(x)=1$ and $\nlc(x)=0$.
As Figure~\ref{fig:sv_jacobian} shows there is indeed an exponential
decay of the singular values and the initial steep decay of the
largest values means that under even relatively small noise in the data
it will be difficult to use more than about ten relevant modes as 
the above reconstruction figures demonstrate.
However the decay rate of the singular values overall is actually more
favourable for reconstructions than for classical exponentially ill-posed problems
such as the backwards or sideways heat problems, the Cauchy problem for
the Laplacian or inverse obstacle scattering.

All of the above reconstructions were obtained using the value
$\alpha=\frac{1}{2}$ \rerevision{ in the fractional damping model}, 
\Margin{Ref 1 3.}
but none were sensitive to this parameter
except for the extreme ends of its range.
However, it is certainly the situation if we had $\alpha=1$ and thus
exponential damping the usefulness of the resulting very small values
of $u$ obtained from anything beyond modest time values would be
extremely limited -- in particular for the reconstruction of $\nlc(x)$
as this is inherently coupled to the magnitude of $u$.
Also a damping term $(\triangle)^\beta$ with $\beta<1$ would weaken the damping and thus influence the degree of ill-posedness. However, theoretically, exponential decay with time persists also with $\beta<1$ as long as $\alpha=1$. We point to \cite{nonlinearity_imaging_fracWest} for an illustration of the pole locations (as being responsible for the degree of ill-posedness) with varying $\beta$ in case of reconstructing $\nlc$ alone.

Two-dimensional reconstructions in the practically relevant case of a piecewise constant coefficient $\nlc$, corresponding to inclusions in a homogeneous background, can be found in \cite{nonlinearity_imaging_2d}.

\section*{Acknowledgment}
The work of the first author was supported by the Austrian Science Fund {\sc fwf} under the grant P36318.
The work of the second author was supported in part by the National Science Foundation through award DMS -2111020.

\end{document}